\documentclass[aap,MSNbibl,dvips]{arximspdf}
\usepackage{graphicx}
\doi{10.1214/13-AAP986}
\volume{25}
\issue{1}
\pubyear{2015}
\firstpage{46}
\lastpage{80}

\makeatletter

\newcommand{\rrvert}{\vert}
\newcommand{\llvert}{\vert}
\newtheorem{theo}{Theorem}
\newtheorem{prop}[theo]{Proposition}
\newtheorem{lem}{Lemma}
\newproclaim{defn}{Definition}
\newproclaim{rem}{Remark}
\newproclaim{assm}{Assumption}

\newcommand{\acal}{{\mathcal A}}
\newcommand{\bcal}{{\mathcal B}}
\newcommand{\ccal}{{\mathcal C}}
\newcommand{\scal}{{\mathcal S}}
\newcommand{\wcal}{{\mathcal W}}

\newcommand{\F}{{\mathcal F}}
\newcommand{\G}{{\mathcal G}}
\newcommand{\Lop}{{\mathcal L}}
\newcommand{\Om}{\Omega}
\newcommand{\R}{{\mathbb R}}
\newcommand{\contr}{\mathrm{c}}
\newcommand{\stopp}{\mathrm{s}}
\newcommand{\e}{{\mathbb{E}}}
\newcommand{\p}{{\mathbb{P}}}
\newcommand{\Q}{{\mathbb{Q}}}

\newcommand{\fraca}[2]{{#1}/{#2}}
\newcommand{\fracf}[2]{({#1})/({#2})}
\makeatother

\begin{document}
\begin{frontmatter}

\title{A zero-sum game between a singular stochastic controller
and a discretionary stopper}
\runtitle{A zero-sum game}

\begin{aug}
\author[A]{\fnms{Daniel} \snm{Hernandez-Hernandez}\ead[label=e1]{dher@cimat.mx}},
\author[B]{\fnms{Robert S.} \snm{Simon}\ead[label=e2]{r.s.simon@lse.ac.uk}}
\and\\
\author[B]{\fnms{Mihail}~\snm{Zervos}\corref{}\ead[label=e3]{mihalis.zervos@gmail.com}}
\runauthor{D. Hernandez-Hernandez, R. S. Simon and
M. Zervos}
\affiliation{Centro de
Investigaci\'{o}n en Matem\'{a}ticas, London School of Economics and
London School of Economics}
\address[A]{D. Hernandez-Hernandez\\
Centro de
Investigaci\'{o}n en Matem\'{a}ticas\\ Apartado Postal 402\\
Guanajuato GTO 36000\\
Mexico\\
\printead{e1}} 
\address[B]{R. S. Simon\\
M. Zervos\\
Department of Mathematics\\
London School of Economics\\ Houghton Street\\
London
WC2A 2AE\\
United Kingdom\\
\printead{e2}\\
\hphantom{\textsc{E-mail:} }\printead*{e3}}
\end{aug}

\received{\smonth{11} \syear{2012}}
\revised{\smonth{6} \syear{2013}}

%
\begin{abstract}
We consider a stochastic differential equation that is controlled
by means of an additive finite-variation process.
A singular stochastic controller, who is a minimizer, determines
this finite-variation process, while a discretionary stopper, who
is a maximizer, chooses a stopping time at which the game
terminates.
We consider two closely related games that are differentiated
by whether the controller or the stopper has a first-move
advantage.
The games' performance indices involve a running payoff as well
as a terminal payoff and penalize control effort expenditure.
We derive a set of variational inequalities that can fully
characterize the games' value functions as well as yield
Markovian optimal strategies.
In particular, we derive the explicit solutions to two special
cases and we show that, in general, the games' value functions
fail to be $C^1$.
The nonuniqueness of the optimal strategy is an interesting
feature of the game in which the controller has the first-move
advantage.
\end{abstract}

%

\end{frontmatter}

\section{Introduction}\label{sec1}

We consider a one-dimensional c\`{a}gl\`{a}d process $X$
that satisfies the stochastic differential equation
%
\begin{equation}
dX_t = b(X_t) \, dt + d\xi_t +
\sigma(X_t) \, dW_t , \qquad X_0 = x \in\R,
\label{X}
\end{equation}
where $\xi$ is a c\`{a}gl\`{a}d finite variation adapted process
such that $\xi_0 = 0$, and $W$ is a standard one-dimensional
Brownian motion.
The games that we analyze involve a controller, who is a
minimizer and chooses a process $\xi$, and a stopper, who
is a maximizer and chooses a stopping time $\tau$.
The two agents share the same performance criterion, which
is given either by
%
\begin{equation}
J_x^v (\xi, \tau) = \e \biggl[ \int
_0^\tau e^{-\Lambda_ t} h(X_t) \, dt
+ \int_{[0,\tau[} e^{-\Lambda_t} \, d\check{\xi}_t +
e^{-\Lambda_\tau} g(X_\tau) \mathbf{1} _{\{ \tau< \infty\}} \biggr]
\label{Jv}
\end{equation}
or by
%
\begin{equation}\quad
J_x^w (\xi, \tau) = \e \biggl[ \int
_0^\tau e^{-\Lambda_ t} h(X_t) \, dt
+ \int_{[0,\tau]} e^{-\Lambda_t} \, d\check{\xi}_t +
e^{-\Lambda_\tau} g(X_{\tau+}) \mathbf{1} _{\{ \tau< \infty\}} \biggr] ,
\label{Jw}
\end{equation}
where $\check{\xi}$ is the total variation process of $\xi$
and
%
\begin{equation}
\Lambda_ t = \int_0^t \bar{\delta}
(X_s) \, ds   \label{Lambda}
\end{equation}
for some positive functions $h, g, \bar{\delta} \dvtx  \R\rightarrow
\R_+$.
The performance index $J^v$ reflects a situation where
the stopper has the ``first-move advantage'' relative to the
controller.
Indeed, if the controller makes a choice such that $\Delta
\xi_0 \neq0$ and the stopper chooses $\tau= 0$, then
$J_x^v (\xi, \tau) = g(x)$.
On the other hand, the performance index $J^w$ reflects
a situation where the controller has the ``first-move advantage''
relative to the stopper: if the controller makes a choice such
that $\Delta\xi_0 \neq0$, and the stopper chooses $\tau= 0$,
then $J_x^w (\xi, \tau) = \llvert \Delta\xi_0\rrvert  + g(x + \Delta\xi_0)$.

Given an initial condition $x \in\R$, $(\xi^*, \tau^*)$ is an
optimal strategy if
%
\begin{equation}
J_x^f \bigl(\xi^*, \tau\bigr) \leq J_x^f
\bigl(\xi^*, \tau^*\bigr) \leq J_x^f \bigl(\xi, \tau^*
\bigr) \label{J*}
\end{equation}
for all admissible strategies $(\xi, \tau)$, where ``$f$''
stands for either ``$v$'' or ``$w$.''
If optimal strategies $(\xi_v^*, \tau_v^*)$, $(\xi_w^*, \tau_w^*)$
exist for the two games for every initial condition
$x \in\R$, then we define the games' value functions by
%
\begin{equation}
v(x) = J_x^v \bigl(\xi_v^*,
\tau_v^*\bigr) \quad\mbox{and} \quad w(x) = J_x^w
\bigl(\xi_w^*, \tau_w^*\bigr) , \label{vw}
\end{equation}
respectively.

Zero-sum games involving a controller and a stopper were
originally studied by Maitra and Sudderth \cite{MS} in a discrete
time setting.
Later, Karatzas and Sudderth \cite{KS} derived the explicit
solution to a game in which the state process is a one-dimensional
diffusion with absorption at the endpoints of a bounded
interval, while, Weerasinghe \cite{W} derived the explicit
solution to a similar game in which the controlled volatility
is allowed to vanish.
Karatzas and Zamfirescu \cite{KZ2} developed a martingale
approach to general controller and stopper games, while,
Bayraktar and Huang \cite{BH} showed that the value function
of such games is the unique viscosity solution to an appropriate
Hamilton--Jacobi--Bellman equation if the state process is a
controlled multi-dimensional diffusion.
Further games involving control as well as discretionary stopping
have been studied by Hamad\`{e}ne and Lepeltier \cite{HL}
and Hamad\`{e}ne \cite{H}.
To a large extent, controller and stopper games have been
motivated by several applications in mathematical finance and
insurance, including the pricing and hedging of American
contingent claims (e.g., see Karatzas and Wang \cite{KW})
and the minimization of the lifetime ruin probability; for example,
see Bayraktar and Young \cite{BY}.

Games such as the ones we study here arise in the context of
several applications.
To fix ideas, consider the singular stochastic control problem
that aims at minimizing the performance criterion
\[
\e \biggl[ \int_0^\infty e^{-\Lambda_ t}
h(X_t) \, dt + \int_{[0,\infty[} e^{-\Lambda_t} \, d
\check{\xi}_t \biggr]
\]
over all controlled processes $\xi$ subject to the dynamics
given by (\ref{X}).
The solution to the special case of this problem that arises when
$b \equiv0$, $\sigma\equiv1$, $\bar{\delta} > 0$ is a constant
and $h(x) = \kappa x^2$, for some $\kappa> 0$, was derived
by Karatzas \cite{K} and is characterized by a constant $\beta$:
it is optimal to exercise minimal control so as to keep the
state process $X$ inside the range $[-\beta, \beta]$ at all times.
The qualitative nature of such a solution has lead to the study
of several applications in which one wants to keep a state process
within an optimal range by means of singular stochastic control.
Such applications include: spaceship control (see Bather
and Chernoff \cite{BC} who introduced singular stochastic
control) where, for example, $X$ represents the deviation of a satellite
from a given altitude and $\xi$ represents fuel expenditure;
the control of an exchange rate (see Miller and Zhang \cite{MZ})
or an inflation rate (see Chiarolla and Haussmann \cite{CH})
where, for example, $X$ models a rate or the fluctuations of a rate
around a target, and $\xi$ models the central bank's cumulative
intervention efforts; the so-called goodwill problem
(see Jack, Jonhnson and Zervos \cite{JJZ}) where, for example,
$X$ is used to model the image that a product has in a
market, and $\xi$ represents the cumulative costs associated
with raising the product's image, for example, through
advertising.\footnote{We have
included here only one indicative reference for each of the
areas mentioned because there is a rich literature
for each of them.}

Any of the applications discussed in the previous paragraph
can give rise to a zero-sum game between a controller and
a stopper that are different incarnations of the same decision
maker.
Such games in which the players model competing
objectives of the same decision maker have attracted
considerable interest in the context of several applications.
For instance, they have been studied in the context
of robust optimization where ``the agent maximizes utility by
his choice of control, while an evil agent minimizes utility by
his choice of perturbation'' (Williams \cite{NW}), or
in the context of time-consistent optimization where
a decision maker's problem is analyzed as a ``game between
successive selves, each of whom can commit for an infinitesimally
small amount of time'' (Ekeland, Mbodji and Pirvu \cite{EMP}).
In what follows, we focus on one of the applications of the
games that we study (several others arising in the context
of the ones discussed in the previous paragraph can be
developed following similar arguments).

Consider a central bank that intervenes to keep fluctuations
of an exchange rate within an optimal range.
At any time, the central bank could be confronted with the
costs of their policy, in particular, with the demand that
its board should be replaced.
In this context, the controller can represent the central bank's
targeting efforts, while the stopper can represent a
political veto on their policy.
In abstract terms, such a problem can be viewed as one of
optimization by a single agent.
However, its analysis and solution requires its formulation
as a zero-sum game.
Indeed, the conflicting natures of such a decision maker's
objectives do not really allow for them to be addressed
by solving a (one-player) stochastic optimization problem.
For instance, the solution to the one-player problem derived by
Davis and Zervos \cite{DZ94}, which is akin to the special
case we solve in Section~\ref{12121212}, involves markedly
different optimal strategies that would be absurd in the context
of an applications such as the one we discuss here.

In particular, the controller tries to minimize, for example, the
performance index $J^w$ given by (\ref{Jw}).
From the controller's perspective, $J^w$ penalizes large
fluctuations of the targeted rate for choices such as
$h(x) = \kappa x^2$, for some $\kappa> 0$, as well as
the expenditure of intervention effort.
On the other hand, the stopper tries to maximize the same
performance criterion $J^w$ because large values of $J^w$
indicate that intervention is ``expensive,'' namely,
unsustainable.
From the stopper's perspective, the choice of the reward
function $g$ can be used to further quantify the bank's reluctance
to intervene, for example, in situations where the rate assumes values
way off the target.
Furthermore, the choice of $J^w$ rather than $J^v$
can be associated with a central bank that is more, rather
than less, keen to intervene.

The development of a theory for zero-sum games such
as the ones we study can therefore provide a useful
analytic tool to decision makers such as a central bank
in their considerations on whether and how to optimally
target a state process such as an exchange rate.
Such analytic tools can be most valuable because getting
a policy wrong can have rather extreme economic
and political consequences.
For instance, one can recall the UK's crash out of the
European Exchange Rate Mechanism (ERM) in 1992.

The games that we study here are the very first ones involving
singular stochastic control and discretionary stopping.
Combining the intuition underlying the solution of standard
singular stochastic control problems and standard optimal
stopping problems by means of variational inequalities (e.g.,
see Karatzas \cite{K} and Peskir and Shiryaev \cite{PS},
resp.), we derive a system of inequalities that can fully
characterize the value function $w$.
We further show that these inequalities can also characterize
the value function $v$ as well as an optimal strategy.
Surprisingly, we have not seen a way to combine all of them
into a single equation.
Our main results include the proof of a verification theorem
that establishes sufficient conditions for a solution to these
inequalities to identify with the value function $w$ and yield
the value function $v$ as well as an optimal strategy, which
we fully characterize.
In this context, we also show that the two games we consider
share the same optimal strategy, and we prove that
\[
v(x) = \max \bigl\{ w(x), g(x) \bigr\} \qquad \mbox{for all } x \in\R.
\]
The nonuniqueness of the optimal strategy when the
controller has the first-move advantage is an interesting
result that arises from our analysis; see Remark~\ref{rem1}
at the end of Section~\ref{vthm-sec}.

We then derive the explicit solutions to two special cases.
The first one is the special case that arises if $X$ is a
standard Brownian motion, and $h$, $g$ are quadratics.
In this case, the value function $w$ is $C^1$, but the $C^1$
regularity of the value function $v$ may fail at a couple
of points.
The second special case is a simpler example revealing
that both of the value functions $w$ and $v$ may fail to be
$C^1$ at certain points and showing that the optimal strategy
may take qualitatively different form, depending on
parameter values.

The paper is organized as follows.
Notation and assumptions are described in Section~\ref{sec2},
while, a heuristic derivation of the system of inequalities
characterizing the solution to the two games is developed
in Section~\ref{sec3} (see Definition~\ref{u}).
In Section~\ref{sec4}, the main results of the paper, namely, a
verification theorem (Theorem~\ref{prop:VT}) and the
construction of the optimal controlled process associated
with a function satisfying the requirements of Definition~\ref{u}
(Lemma~\ref{lem:xi*tau*}) are proved.
In Sections~\ref{12121212} and~\ref{16161616}, the explicit
solutions to two nontrivial special cases are derived.

\section{Notation and assumptions}\label{sec2}
\label{assm}

We fix a filtered probability space $(\Om, \F, \allowbreak \F_t , \p)$ satisfying
the usual conditions and carrying a standard one-dimensional
$(\F_t)$-Brownian motion $W$.
We denote by $\acal_\stopp$ the set of all $(\F_t)$-stopping
times and by $\acal_\contr$ the family of all $(\F_t)$-adapted
finite-variation c\`{a}gl\`{a}d processes $\xi$ such that $\xi_0 = 0$.
Every process $\xi\in\acal_\contr$ admits the decomposition
$\xi= \xi^{\mathrm c} + \xi^{\mathrm j}$ where
$\xi^{\mathrm c}$, $ \xi^{\mathrm j}$ are $(\F_t)$-adapted
finite-variation c\`{a}gl\`{a}d processes such that
$\xi^{\mathrm c}$ has continuous sample paths,
\[
\xi_0^{\mathrm c} = \xi_0^{\mathrm j} = 0 \quad
\mbox{and} \quad\xi_t^{\mathrm j} = \sum
_{0 \leq s < t} \Delta\xi_s \qquad \mbox{for all } t>0 ,
\]
where $\Delta\xi_s = \xi_{s+} - \xi_s$ for $s \geq0$.
Given such a decomposition, there exist $(\F_t)$-adapted
continuous processes $(\xi^{\mathrm c})^+$,
$(\xi^{\mathrm c})^-$ such that
\[
\bigl(\xi^{\mathrm c}\bigr)_0^+ = \bigl(\xi^{\mathrm c}
\bigr)_0^- =0 , \qquad \xi^{\mathrm c} = \bigl(\xi^{\mathrm c}
\bigr)^+ - \bigl(\xi^{\mathrm c}\bigr)^- \quad\mbox{and} \quad \check{
\xi}^{\mathrm c} = \bigl(\xi^{\mathrm c}\bigr)^+ + \bigl(\xi^{\mathrm c}
\bigr)^- ,
\]
where $\check{\xi}^{\mathrm c}$ is the total variation
process of $\xi^{\mathrm c}$.

The following assumption that we make implies that, given
any $\xi\in\acal_\contr$, (\ref{X}) has a unique strong solution;
see Protter \cite{P}, Theorem~V.7.

\begin{assm} \label{Assm1}
The functions $b, \sigma\dvtx  \R\rightarrow\R$
satisfy
\[
\bigl\llvert b(x) - b(y)\bigr\rrvert + \bigl\llvert \sigma(x) - \sigma(y)\bigr
\rrvert \leq K \llvert x-y\rrvert \qquad \mbox{for all } x,y \in\R,
\]
for some constant $K>0$, and $\sigma^2 (x) > \sigma_0$
for all $x \in\R$, for some constant $\sigma_0 > 0$.
\end{assm}

We also make the following assumption on the data of the
reward functionals defined by (\ref{Jv})--(\ref{Lambda}).

\begin{assm} \label{Assm2}
The functions $\bar{\delta} , h , g\dvtx  \R\rightarrow\R_+$ are
continuous, and there exists a constant $\delta> 0$
such that $\bar{\delta} (x) > \delta$ for all $x \in\R$.
\end{assm}

It is worth noting at this point that, given $\xi\in\acal_\contr$,
we may have $\e[\check{\xi}_t] = \infty$, for some $t>0$.
In such a case, the reward functionals given by
(\ref{Jv})--(\ref{Jw}) are well defined but may take the
value $\infty$.

\section{Heuristic derivation of variational inequalities for the value function \texorpdfstring{$\bolds{w}$}{w}}\label{sec3}
\label{heuristics}

Before addressing the game, we consider the optimization
problems faced by the two players in the absence of
competition.
To this end, we consider any bounded interval $]\gamma_1,
\gamma_2[$, we denote by $T_{\gamma_1}$ (resp.,
$T_{\gamma_2}$) the first hitting time of $\{ \gamma_1 \}$
(resp., $\{ \gamma_2 \}$), and we fix any constants
$C_{\gamma_1}, C_{\gamma_2} \geq0$.

Given an initial condition $x \in\, ]\gamma_1,
\gamma_2 [$, a controller is concerned with solving the singular
stochastic control problem whose value function is given by
%
\begin{eqnarray}\label{u-ssc}
&&v _{\mathrm{ssc}} (x; \gamma_1, \gamma_2,
C_{\gamma_1}, C_{\gamma_2})\nonumber\\
 &&\qquad = \inf_{\xi\in\acal_\contr} \e \biggl[ \int
_0^{T_{\gamma_1} \wedge T_{\gamma_2}} e^{-\Lambda_t} h(X_t) \, dt
+ \int_{[0, T_{\gamma_1} \wedge T_{\gamma_2}[} e^{-\Lambda_t} \, d\check{\xi}_t
\\
&&\qquad \quad \hphantom{\inf_{\xi\in\acal_\contr} \e \biggl[} {} + e^{-\Lambda_{T_{\gamma_1}}}
C_{\gamma_1} \mathbf{1} _{\{ T_{\gamma_1} < T_{\gamma_2} \}} + e^{-\Lambda
_{T_{\gamma_2}}} C_{\gamma_2}
\mathbf{1} _{\{ T_{\gamma_2}
< T_{\gamma_1} \}} \biggr] .
\nonumber
\end{eqnarray}
In the presence of Assumptions~\ref{Assm1} and~\ref{Assm2},
$v _{\mathrm{ssc}}$ is $C^1$ with absolutely continuous
first derivative and identifies with the solution to the
variational inequality
\[
\min \bigl\{ \Lop u(x) + h(x) ,  1 - \bigl\llvert u'(x)\bigr
\rrvert \bigr\} = 0
\]
with boundary conditions
\[
u(\gamma_1) = C_{\gamma_1} \quad\mbox{and} \quad u(
\gamma_2) = C_{\gamma_2} ,
\]
where the operator $\Lop$ is defined by
%
\begin{equation}
\Lop u(x) = {\tfrac{1}{2}}\sigma^2 (x) u''(x)
+ b(x) u'(x) - \bar {\delta} (x) u(x) \label{L};
\end{equation}
see Sun \cite{S}, Theorem~3.2.
In this case, it is optimal to exercise minimal action so
that the state process $X$ is kept outside the interior
of the set
\[
\ccal_{\mathrm{ssc}} = \bigl\{ x \in\, ]\gamma_1,
\gamma_2[\, \mid \bigl\llvert u'(x)\bigr\rrvert =
1 \bigr\} .
\]

Given an initial condition $x \in\, ]\gamma_1,
\gamma_2[$, a stopper faces the discretionary stopping
problem whose value function is given by
%
\begin{eqnarray}\label{u-ds}
&&v _{\mathrm{ds}} (x; \gamma_1, \gamma_2,
C_{\gamma_1}, C_{\gamma_2})\nonumber\\
 &&\qquad = \sup_{\tau\in\acal_\stopp} \e \biggl[ \int
_0^{\tau\wedge T_{\gamma_1} \wedge T_{\gamma_2}} e^{-\Lambda_t} h(X_t) \, dt
+ e^{-\Lambda_\tau} g(X_\tau) \mathbf{1} _{\{ \tau< T_{\gamma_1} \wedge T_{\gamma_2} \}}
\\
&&\qquad \quad \hphantom{\sup_{\tau\in\acal_\stopp} \e \biggl[\;\;\;} {} + e^{-\Lambda_{T_{\gamma_1}}}
C_{\gamma_1} \mathbf{1} _{\{ T_{\gamma_1} \leq\tau\wedge T_{\gamma_2} \}} + e^{-\Lambda_{T_{\gamma_2}}} C_{\gamma_2}
\mathbf{1} _{\{ T_{\gamma_2} \leq\tau\wedge T_{\gamma_1} \}} \biggr] ,
\nonumber
\end{eqnarray}
where $X$ is the solution to (\ref{X}) for $\xi\equiv0$.
In this case, Assumptions~\ref{Assm1} and~\ref{Assm2}
ensure that $v_{\mathrm{ds}}$ is the difference of two
convex functions and identifies with the solution, in an
appropriate distributional sense, to the variational inequality
\[
\max \bigl\{ \Lop u(x) + h(x) ,  g(x) - u(x) \bigr\} = 0
\]
with boundary conditions
\[
u(\gamma_1) = C_{\gamma_1} \quad\mbox{and} \quad u(
\gamma_2) = C_{\gamma_2} ,
\]
where $\Lop$ is defined by (\ref{L}); see Lamberton and
Zervos \cite{LZ}, Theorems~12 and 13.
In this case, the optimal stopping time $\tau^\circ$
identifies with the first hitting time of the so-called stopping
region
\[
{\mathcal S} _{\mathrm{ds}} = \bigl\{ x \in\, ]\gamma_1,
\gamma_2[ \, \mid u(x) = g(x) \bigr\} ,
\]
namely, $\tau^\circ= \inf\{ t \geq0 \mid X_t \in
{\mathcal S} _{\mathrm{ds}} \}$.

Now, we consider the game where the controller has
the ``first-move advantage'' relative to the stopper, and we
assume that there exists a Markovian optimal
strategy $(\xi^*, \tau^*)$ for the sake of the discussion
in this section.
We expect that this optimal strategy involves the
same tactics as the ones we have discussed above.
From the perspective of the controller, the state space
$\R$ splits into a control region $\mathcal C$ and a
waiting region ${\mathcal W} _{\mathrm c}$.
Accordingly, $\xi^*$ should involve minimal action to
keep the state process in the closure $\R\setminus
\operatorname{int} {\mathcal C}$ of the waiting region
${\mathcal W} _{\mathrm c}$ for as long as the stopper
does not terminate the game.
Similarly, from the perspective of the stopper, the state
space $\R$ splits into a stopping region $\mathcal S$
and a waiting region ${\mathcal W} _{\mathrm s}$,
and $\tau^*$ is the first hitting time of $\mathcal S$.

Inside any bounded interval $]\gamma_1, \gamma_2[
\, \subseteq{\mathcal W} _{\mathrm s}$, the
requirement that $(\xi^*, \tau^*)$ should satisfy (\ref{J*})
suggests that $w$ should identify with $v _{\mathrm{ssc}}$
defined by (\ref{u-ssc}) for $C_{\gamma_1} = w(\gamma_1)$
and $C_{\gamma_2} = w(\gamma_2)$.
Therefore, we expect that $w$ should satisfy
%
\begin{equation}
\min \bigl\{ \Lop w(x) + h(x) ,  1 - \bigl\llvert w'(x)\bigr
\rrvert \bigr\} = 0 \qquad \mbox{inside } {\mathcal W} _{\mathrm s} .
\label{VIh1}
\end{equation}
Inside any bounded interval $]\gamma_1, \gamma_2[
\, \subseteq{\mathcal W} _{\mathrm c}$, the
requirement that $(\xi^*, \tau^*)$ should satisfy (\ref{J*})
suggests that $w$ should identify with $v _{\mathrm{ds}}$
defined by (\ref{u-ds}) for $C_{\gamma_1} = w(\gamma_1)$
and $C_{\gamma_2} = w(\gamma_2)$.
Therefore, we expect that $w$ should satisfy
%
\begin{equation}
\max \bigl\{ \Lop w(x) + h(x) ,  g(x) - w(x) \bigr\} = 0 \qquad \mbox{inside }
{\mathcal W} _{\mathrm c} . \label{VIh2}
\end{equation}

To couple variational inequalities (\ref{VIh1}) and
(\ref{VIh2}), we consider four possibilities.
The region $\wcal\wcal= \wcal_{\mathrm c} \cap
\wcal_{\mathrm s}$ where both players should
wait is associated with the inequalities
%
\begin{equation}
\Lop w + h = 0 , \qquad\bigl\llvert w'\bigr\rrvert < 1 \quad
\mbox{and} \quad g < w . \label{v-ineq1}
\end{equation}
Inside the set $\ccal\wcal= \ccal\cap\wcal_{\mathrm s}$
where the stopper should wait, whereas,
the controller should act, we expect that
%
\begin{equation}
\Lop w + h \geq0 , \qquad\bigl\llvert w'\bigr\rrvert = 1 \quad
\mbox{and} \quad g < w .
\end{equation}
Inside the part of the state space $\wcal\scal=
\wcal_{\mathrm c} \cap\scal$ where
the controller would rather wait if the stopper deviated
from the optimal strategy and did not terminate the
game, we expect that
%
\begin{equation}
\Lop w + h \leq0 , \qquad\bigl\llvert w'\bigr\rrvert < 1 \quad
\mbox{and} \quad g = w .
\end{equation}
Finally, the region $\ccal\scal= \ccal\cap\scal$ in which
the stopper should terminate the game should the controller
deviate from the optimal strategy and did not act, we expect that
%
\begin{equation}
\Lop w + h \in\R, \qquad\bigl\llvert w'\bigr\rrvert = 1 \quad
\mbox{and} \quad g \geq w .
\end{equation}
These inequalities give rise to the following definition.
Here, as well as in the rest of the paper, we denote by
$\operatorname{int} \Gamma$ and $\operatorname{cl}
\Gamma$ the interior and the closure of a set $\Gamma
\subseteq\R$, respectively.

%
\begin{defn} \label{u}
A \emph{candidate for the value function} $w$ is a continuous
function $u\dvtx  \R\rightarrow\R_+$ that is $C^1$ with absolutely
continuous first derivative inside $\R\setminus{\mathcal B}$,
where ${\mathcal B}$ is a finite set, satisfies
\[
\bigl\llvert u'(x)\bigr\rrvert \leq1 \qquad\mbox{for all } x \in\R
\setminus {\mathcal B} ,
\]
and has the following properties, where
\begin{eqnarray*}
\ccal&=& \operatorname{cl} \bigl[ \operatorname{int} \bigl\{ x \in\R\setminus
\bcal\mid \bigl\llvert u'(x)\bigr\rrvert = 1 \bigr\} \bigr] ,
\\
\scal_\wcal&=& \bigl\{ x \in\R\mid u(x) = g(x) \bigr\} , \qquad
\scal_\ccal= \operatorname{cl} \bigl\{ x \in \R\mid u(x) < g(x) \bigr\}
,
\\
\scal&=& \scal_\wcal\cup\scal_\ccal\quad\mbox{and} \quad
\wcal= \R\setminus ( \ccal\cup\scal ) .
\end{eqnarray*}

\begin{longlist}[(III)]
\item[(I)] Each of the sets $\ccal$, $\scal_\wcal$ and
$\scal_\ccal$ is a finite union of intervals, and
$\bcal\subseteq\allowbreak \scal_\ccal\subseteq\ccal$.

\item[(II)] $u$ satisfies
\[
\Lop u(x) + h(x) %
\cases{\displaystyle = 0 ,&\quad Lebesgue-a.e. in $
\wcal$,
\cr
\displaystyle \geq0 , &\quad Lebesgue-a.e. in $\operatorname{int} (
\ccal\setminus\scal ) $,
\cr
\displaystyle \equiv\Lop g(x) + h(x) \leq0 ,&\quad
Lebesgue-a.e. in $\operatorname{int} \scal_\wcal$,
\cr
\displaystyle \in
\R, &\quad Lebesgue-a.e. in $\operatorname{int} \scal_\ccal\setminus
\bcal$. } %
\]

\item[(III)] If we denote by $u_-' (c)$ [resp., $u_+' (c)$] the left-hand
(resp., right-hand) derivative of $u$ at $c \in{\mathcal B}$,
then
\[
\mbox{either} \qquad u_-' (c) = 1 \quad \mbox{and}\quad
u_+' (c) < 1 \quad\mbox{or} \quad u_-' (c) > -1 \quad
\mbox{and}\quad u_+' (c) = -1
\]
for all $c \in{\mathcal B}$.
\end{longlist}
\end{defn}

In the following definition, we introduce some terminology
we are going to use.

%
\begin{defn} \label{terminology}
Given a function $u$ satisfying the conditions of Definition~\ref{u},
we call the regions ${\mathcal W}$, ${\mathcal C}$ and ${\mathcal S}$
\emph{waiting}, \emph{control} and \emph{stopping}, respectively.
Also, we call \emph{reflecting} all finite boundary points $x$ of
${\mathcal C}$ such that
\[
u' (x-\varepsilon) < 1 \quad \mbox{and}\quad  u'(x) = 1 \quad
\mbox{or} \quad u' (x) = -1 \quad \mbox{and}\quad  u'(x+
\varepsilon) > -1
\]
for all $\varepsilon> 0$ sufficiently small, and \emph{repelling}
all other finite boundary points of ${\mathcal C}$.
\end{defn}

It is worth noting that requirement (III) of Definition~\ref{u}
implies that all points in ${\mathcal B}$ are repelling.
The special case that we solve in Section~\ref{12121212}
involves only reflecting boundary points.
On the other hand, the special case that we solve in
Section~\ref{16161616} involves repelling as well as
reflecting points and ${\mathcal B} \neq\varnothing$.

\section{A verification theorem}\label{sec4}
\label{vthm-sec}

Before addressing the main result on this section, namely
Theorem~\ref{prop:VT}, we consider the following result,
which is concerned with the construction of the process
$\xi^*$ that is part of the optimal strategy associated with
a given function satisfying the requirements of
Definition~\ref{u}.
The main idea of its proof is to paste solutions to (\ref{X})
that are reflecting in appropriate boundary points.

%
\begin{lem} \label{lem:xi*tau*}
Consider a function $u\dvtx  \R\rightarrow\R_+$ that satisfies
the conditions of Definition~\ref{u}.
There exists a controlled process $\xi^* \in\acal_\contr$
such that
%
\begin{eqnarray}\label{xi*1}
&\displaystyle \mbox{the set } \bigl\{ t \geq0 \mid X_t^* \in{
\mathcal B} \bigr\} \mbox{ is finite} ,&
\\\label{xi*2}
&\displaystyle X_t^* \in\R\setminus\operatorname{int} {\mathcal C}
\qquad \mbox{for all } t > 0 , \qquad u\bigl(X_{t+}^*\bigr) - u
\bigl(X_t^*\bigr) = - \bigl\llvert \Delta\xi_t^*\bigr
\rrvert = - \bigl\llvert \Delta X_t^*\bigr\rrvert& \nonumber
\\[-4pt]
\\[-12pt]
\eqntext{\mbox{for all
} t \geq0 ,}
\\\label{xi*3}
&\displaystyle \bigl(\xi^{*\mathrm{c}}\bigr)_t^+ = \int
_0^t \mathbf{1} _{\{ u' (X_s^*)
= -1 \}} \, d\bigl(
\xi^{*\mathrm{c}}\bigr)_s^+ \quad \mbox{and} \quad \bigl(
\xi^{*\mathrm{c}}\bigr)_t^- = \int_0^t
\mathbf{1} _{\{ u' (X_s^*) =
1 \}} \, d\bigl(\xi^{*\mathrm{c}}\bigr)_s^-&
\nonumber
\\[-6pt]
\\[-10pt]
\eqntext{\mbox{for all } t \geq0 , }
\end{eqnarray}
where $X^*$ is the associated solution to (\ref{X}).
\end{lem}

\begin{pf}
Given a finite interval $[\alpha, \beta]$ and a controlled process
$\xi\in\acal_\contr$, suppose that there exists a point $\bar{x}
\in[\alpha, \beta]$ and an $(\F_t)$-stopping time $\tau$ with
$\p(\tau< \infty) > 0$ such that the solution to (\ref{X}) is such
that $X_\tau= \bar{x}$ on\vadjust{\goodbreak} the event $\{ \tau< \infty\}$.
On the probability space $(\Om, \F, \G_t, \Q)$, where $(\G_t)$
is the filtration defined by $\G_t = \F_{\tau+t}$ and $\Q$ is
the conditional probability measure $\p(\cdot\mid\tau< \infty)$
that has Radon--Nikodym derivative with respect to $\p$ given by
\[
\frac{d \Q}{d \p} = \frac{1}{\p( \tau< \infty)} \mathbf{1} _{\{
\tau
< \infty\}} ,
\]
the process $B$ defined by $B_t = (W_{\tau+t} - W_\tau)
\mathbf{1} _{\{ \tau< \infty\}}$ is a standard $(\G_t)$-Brownian
motion that is independent of $\G_0 = \F_\tau$; see Revuz
and Yor \cite{RY}, Exercise IV.3.21.
In this context, there exist $(\G_t)$-adapted continuous
processes $\bar{X}$ and $\bar{\xi}$ such that $\bar{\xi}$ is a
finite variation process,
\begin{eqnarray}
\bar{X}_t = \bar{x} + \int_0^t
b(\bar{X}_s) \, ds + \bar{\xi}_t + \int
_0^t \sigma(\bar{X}_s) \,
dB_s ,\nonumber
\\
\eqntext{\displaystyle \bar{X}_t \in[\alpha, \beta] ,   \bar{\xi}_t^+  =
\int_0^t \mathbf{1} _{\{ \bar{X}_s = \alpha\}} \, d\bar{
\xi}_s^+  \mbox{ and }   \bar{\xi}_t^- = \int
_0^t \mathbf{1} _{\{ \bar{X}_s = \beta\}} \, d\bar{
\xi}_s^- ;}
\end{eqnarray}
see El Karoui and Chaleyat-Maurel \cite{EC} and
Schmidt \cite{Sc}.
Since $(t-\tau)^+$ is an $(\F_{\tau+t})$-stopping time,
$\G_t = \F_{\tau+t}$ and $B_{(t-\tau)^+} = (W_t - W_\tau)
\mathbf{1} _{\{ \tau< t \}}$ for all $t \geq0$,
\begin{eqnarray*}
\bar{X} _{(t-\tau)^+} & =& \bar{x} + \int_0^{(t-\tau)^+}
b ( \bar{X} _s ) \, ds + \bar{\xi}_{(t-\tau)^+} + \int
_0^{(t-\tau)^+} \sigma ( \bar{X} _s ) \,
dB_s
\\
& = &\bar{x} + \int_0^t b ( \bar{X}
_{(s-\tau)^+} ) \, d{(s-\tau)^+} + \bar{\xi}_{(t-\tau)^+} + \int
_0^t \sigma ( \bar{X} _{(s-\tau)^+} ) \,
dB_{(s-\tau)^+}
\\
& = &\bar{x} + \int_0^t \mathbf{1}
_{\{ \tau\leq s \}} b ( \bar{X} _{(s-\tau)^+} ) \, ds + \bar{\xi}_{(t-\tau)^+} +
\int_0^t \mathbf{1} _{\{ \tau\leq s \}} \sigma (
\bar{X} _{(s-\tau)^+} ) \, dW_s ;
\end{eqnarray*}
see Revuz and Yor \cite{RY}, Propositions~V.1.4, V.1.5.
Similarly we can see, for example, that
\[
\bar{\xi}_{(t-\tau)^+}^+ = \int_0^{(t-\tau)^+}
\mathbf{1} _{\{
\bar{X}_s
= \alpha\}} \, d\bar{\xi}_s^+ = \int
_0^t \mathbf{1} _{\{ \bar{X}
_{(s-\tau)^+} = \alpha\}} \, d\bar{
\xi}_{(s-\tau)^+}^+ .
\]
In view of this observation, we can see that, if we define
%
\begin{equation}
\tilde{X}_t = %
\cases{\displaystyle X_t ,&
\quad if $t \leq\tau$,
\cr
\displaystyle \bar{X}_{t - \tau} ,&\quad if $t >
\tau$, } %
\quad\mbox{and} \quad \tilde{\xi}_t = %
\cases{\displaystyle \xi_t ,&\quad if $t \leq\tau$,
\cr
\displaystyle
\bar{\xi}_{t - \tau} ,&\quad if $t > \tau$, } %
\label{Rlr1}
\end{equation}
then $\tilde{X}$ is the solution to (\ref{X}) that is driven by
$\tilde{\xi} \in\acal_\contr$,
%
\begin{eqnarray}\label{Rlr2}
\tilde{X}_t \in[\alpha, \beta] , \qquad \tilde{\xi}_t^+
- \tilde{\xi}_{\tau+}^+ &=& \int_\tau^t
\mathbf{1} _{\{ \tilde{X}_s = \alpha\}} \, d\tilde{\xi}_s^+ \quad\mbox{and} \nonumber
\\[-8pt]
\\[-8pt]
\tilde{\xi}_t^- - \tilde{\xi} _{\tau+}^- &=& \int
_\tau^t \mathbf{1} _{\{ \tilde{X}_s = \beta\}} \, d\tilde{
\xi}_s^-
\nonumber
\end{eqnarray}
for all $t > \tau$.\vadjust{\goodbreak}

Using the same arguments and references, we can show that,
given an interval $[\alpha, \infty[$, a point $\bar{x} \in[\alpha,
\infty[$, a controlled process $\xi\in\acal_\contr$ and an
$(\F_t)$-stopping time $\tau$ such that the solution to (\ref{X})
is such that $X_\tau= \bar{x}$ on the event $\{ \tau< \infty\}$,
there exist processes $\tilde{X}$ and $\tilde{\xi} \in\acal_\contr$
satisfying (\ref{X}) and such that
%
\begin{eqnarray}
\label{Rl1} \tilde{X}_t = X_t \quad\mbox{and} \quad
\tilde{\xi} _t = \xi_t \qquad \mbox{for all } t \leq\tau,
\\
\label{Rl2}
\\
[-18.5pt] \eqntext {\displaystyle \tilde{X}_t \in [\alpha, \infty[ , \tilde{
\xi}_t^+ - \tilde{\xi}_{\tau+}^+ = \int_\tau^t
\mathbf{1} _{\{ \tilde{X}_s = \alpha\}} \, d\tilde{\xi}_s^+ \mbox{ and } \tilde{
\xi}_t^- - \tilde{\xi} _{\tau+}^- = 0 \mbox{ for all } t >
\tau.}
\end{eqnarray}
Similarly, given an interval $ ] {-}\infty, \beta ]$, a
point $\bar{x} \in\, ] {-}\infty, \beta ]$, a controlled
process $\xi\in\acal_\contr$ and an $(\F_t)$-stopping time
$\tau$ such that the solution to (\ref{X}) is such that
$X_\tau= \bar{x}$ on the event $\{ \tau< \infty\}$, there
exist processes $\tilde{X}$ and $\tilde{\xi} \in\acal_\contr$
satisfying (\ref{X}) and such that
%
\begin{eqnarray}
\tilde{X}_t = X_t \quad\mbox{and} \quad\tilde{\xi}
_t = \xi_t \qquad \mbox{for all } t \leq\tau,
\label{Rr1}
\\
\label{Rr2}
\\
[-18.5pt] \eqntext{\displaystyle \tilde{X}_t \in \,] {-}\infty, \beta ] , \tilde{
\xi}_t^+ - \tilde{\xi}_{\tau+}^+ = 0 \mbox{ and } \tilde{
\xi}_t^- - \tilde{\xi} _{\tau+}^- = \int_\tau^t
\mathbf{1} _{\{ \tilde{X}_s = \beta\}} \, d\tilde{\xi}_s^- \mbox{ for all } t >
\tau.}
\end{eqnarray}

%
\begin{figure}[b]

\includegraphics{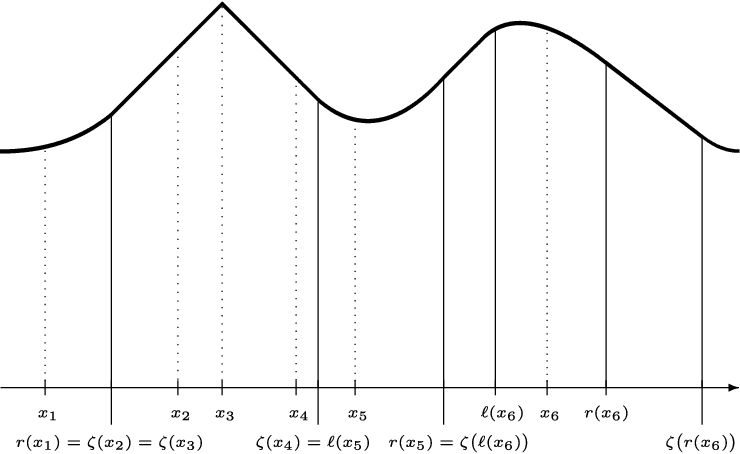}

\caption{Illustration of
the functions $\zeta$, $\ell$, $r$ appearing
in the proof of Lemma~\protect\ref{lem:xi*tau*}. The vertical
solid lines also demarcate the region $\ccal$.}
\label{fig1}\end{figure}

Given a function $u$ that satisfies the requirements of
Definition~\ref{u}, we now use the notation and the
terminology introduced by Definitions~\ref{u}
and~\ref{terminology} to iteratively construct a process
$\xi^* \in\acal_\contr$ such that (\ref{xi*1})--(\ref{xi*3})
hold true by means of the constructions above.
To this end, we introduce the following notation,
which is illustrated by Figure~\ref{fig1}.
If $\operatorname{int} {\mathcal C} \neq\varnothing$ and
$x \in{\mathcal C}$, then we\vadjust{\goodbreak} recall that we use $u_-' (x)$
[resp., $u_+' (x)$] to denote the left-hand (resp., the
right-hand) first derivative of $u$ at $x$, we define
\[
\zeta(x) = %
\cases{\displaystyle \sup \{ y < x \mid y \notin {
\mathcal C} \} ,&\quad if $u_-'(x) = 1 $,
\cr
\displaystyle \inf \{ y
> x \mid y \notin{\mathcal C} \} ,&\quad if $u_+' (x) = -1$ and
$u_-' (x) < 1$, } %
\]
and we note that $\zeta(x) \in\R$ because $u$ is
real-valued.
On the other hand, given any $x \in\R$, we define
\[
\ell(x) = \sup \{ y < x \mid y \in\operatorname{int} {\mathcal C} \} \quad
\mbox{and} \quad r(x) = \inf \{ y > x \mid y \in\operatorname{int} {\mathcal C}
\} ,
\]
with the usual conventions that $\sup\varnothing=
-\infty$ and $\inf\varnothing= \infty$.
The algorithm that we now develop terminates after finite
iterations because each of the sets ${\mathcal C}$,
${\mathcal W}$ is a finite union of intervals.

\emph{STEP 0: Initialization}.
We consider the following four possibilities that can happen,
depending on the initial condition $x$ of (\ref{X}):

If $\operatorname{int} {\mathcal C} \neq\varnothing$ and
$x \in\operatorname{int} {\mathcal C}$ (e.g., see the
points $x_2$, $x_3$, $x_4$ in Figure~\ref{fig1}), then we define
$\xi_t^0 = \zeta(x) - x$ for all $t > 0$.
If we denote by $X^0$ the corresponding solution to
(\ref{X}), and we set $\tau_0 = 0$, then $X^0$ has a single
jump at time $\tau_0$,
\begin{eqnarray*}
\begin{tabular}{@{}l@{}}
 $u(X_{0+}^0) - u(X_0^0
)  =  u ( \zeta(x) ) - u(x) = - \llvert \zeta(x) - x
\rrvert = - \llvert \Delta\xi_0^0\rrvert $,
\\
 $\mbox{if } \zeta(x) < x , \mbox{ then } X_{\tau_0+}^0 =
X_{0+}^0 = \zeta(x) = r ( \zeta(x) ) = r
(X_{0+}^0) \mbox{ is reflecting}$
\\
 $\mbox{and if } x < \zeta(x) , \mbox{ then } X_{\tau_0+}^0 =
X_{0+}^0 = \zeta(x) = \ell ( \zeta(x) ) = \ell
(X_{0+}^0) \mbox{ is reflecting} $.
\end{tabular}
\end{eqnarray*}
In this case, $X_0 \in{\mathcal B}$ if $x \in{\mathcal B}
\subseteq{\mathcal C}$.

If $\ell(x) = -\infty$ and $r(x) = \infty$, which is the
case if $\operatorname{int} {\mathcal C} = \varnothing$,
then we define $\xi^0 = 0$, we denote by $X^0$ the
corresponding solution to (\ref{X}), and we let $\tau_0
= \infty$.

If $\operatorname{int} {\mathcal C} \neq\varnothing$,
$x \in\R\setminus\operatorname{int} {\mathcal C}$
and either of $\ell(x)$, $r(x)$ is reflecting (e.g., see
the points $x_1$, $x_5$ in Figure~\ref{fig1}), then we
define $\xi^0 = 0$, we denote by $X^0$ the
corresponding solution to (\ref{X}), and we set $\tau_0
= 0$.

If $\operatorname{int} {\mathcal C} \neq\varnothing$,
$x \in\R\setminus\operatorname{int} {\mathcal C}$,
and both $\ell(x)$, $r(x)$ are repelling if finite
(e.g., see the point $x_6$ in Figure~\ref{fig1}), then we
consider the $(\F_t)$-stopping times
\[
T_{\ell(x)} = \inf \bigl\{ t \geq0 \mid X_t^\dagger\leq
\ell(x) \bigr\} , \qquad T_{r(x)} = \inf \bigl\{ t \geq0 \mid X_t^\dagger\geq r(x) \bigr\} ,
\]
where $X^\dagger$ is the solution to (\ref{X}) for $\xi= 0$,
and we set
\[
\xi_t^0 = \bigl[ \zeta \bigl( \ell(x) \bigr) - \ell(x)
\bigr] \mathbf{1} _{\{ T_{\ell(x)} < T_{r(x)} \wedge t \}} + \bigl[ \zeta \bigl( r(x) \bigr) - r(x) \bigr]
\mathbf{1} _{\{ T_{r(x)}
< T_{\ell(x)} \wedge t \}} ,
\]
in which expression, we define $\zeta ( \ell(x)  )
- \ell(x)$ [resp., $\zeta ( r(x)  ) - r(x)$] arbitrarily
if $\ell(x) = -\infty$ [resp., $r(x) = \infty$].
If we denote by $X^0$ the corresponding solution to
(\ref{X}), and we set $\tau_0 = T_{r(x)} \wedge T_{\ell(x)}$,
then $X^0$ has a single jump at the $(\F_t)$-stopping
time $\tau_0$,
\begin{eqnarray*}
&\displaystyle  X_t^0 \in\R\setminus\operatorname{int} {\mathcal C}
\quad\mbox{and} \quad u\bigl(X_{t+}^0\bigr) - u
\bigl(X_t^0\bigr) = - \bigl\llvert \Delta
\xi_t^0\bigr\rrvert \qquad \mbox{for all } t \leq
\tau_0 ,&
\\
&\displaystyle \begin{tabular}{@{}l@{}}
 $\mbox{on the event } \{ T_{\ell(x)} < T_{r(x)} \} \in\F
_{\tau_0}, \mbox{ the point } X_{\tau_0+}^0 = \zeta
( \ell (x) ) \mbox{ is reflecting}$
\\
 $\mbox{and on the event } \{ T_{r(x)} < T_{\ell(x)} \} \in\F
_{\tau_0} , \mbox{ the point } X_{\tau_0+}^0 = \zeta
( r(x) )$ is reflec-\\
ting.
\end{tabular}&
\end{eqnarray*}
In this case, we may have $X_{\tau_0}^0 \in{\mathcal B}$
but $X_{\tau_0+}^0 \notin{\mathcal B}$ and $X_t^0 \notin
{\mathcal B}$ for all $t < \tau_0$.\vadjust{\goodbreak}

\emph{STEP 1: Induction hypothesis}.
We assume that we have determined an $(\F_t)$-stopping
time $\tau_j$, and we have constructed a process $\xi^j \in
{\mathcal A} _\mathrm{c}$ such that, if we denote by $X^j$ the
associated solution to (\ref{X}), then (\ref{xi*1})--(\ref{xi*3})
are satisfied for $\xi^j$, $X^j$ in place of $\xi^*$, $X^*$
and for all $t \leq\tau_j$ instead of all positive $t$.
Also, we assume that, if $\p(\tau_j < \infty) > 0$,
then one of the following two possibilities occur:
\begin{longlist}[(II)]
\item[(I)] there exists a point $x^j$ such that $X_{\tau_j}^j =
x^j$ on the event $\{ \tau_j < \infty\}$;

\item[(II)] there exist points $x_1^j$, $x_2^j \in\R$ and events
$A_1^j, A_2^j \in\F_{\tau_j}$\vspace*{-1pt} forming a partition of
$\{ \tau_j < \infty\}$ such that $\p(A_k^j) > 0$,
$X_{\tau_j+}^j = x_k^j$ on the event\vspace*{-1pt} $A_k^j$ and at
least one of $\ell(x_k^j)$, $r(x_k^j)$ is finite and
reflecting, for $k = 1,2$.
\end{longlist}

Step 0 provides such a construction for $j=0$.
In particular, the last possibility there gives rise to
Case~(II) for
\begin{eqnarray*}
A_1^0 &=& \{ T_{\ell(x)} < T_{r(x)} \} ,
\qquad A_2^0 = \{ T_{r(x)} < T_{\ell(x)} \}
,\\
 x_1^0 &=& \zeta \bigl( \ell(x) \bigr) \quad
\mbox{and} \quad x_2^0 = \zeta \bigl( r (x) \bigr) .
\end{eqnarray*}
On the other hand, the second possibility there is
such that $\p(\tau_j < \infty) = 0$, while the
remaining two possibilities give rise to Case~(I).

\emph{STEP 2}.
If $\p(\tau_j < \infty) = 0$, then define $\xi^* = \xi^j$,
$X^* = X^j$ and stop.
Otherwise, we proceed to the next step.

\emph{STEP 3}.
We address the situation arising in the context
of Case~(II) of Step~1; the analysis regarding Case~(I)
is simpler and follows exactly the same steps.
To this end, we first consider the $(\F_t)$-stopping
time $\hat{\tau} = \tau_j \mathbf{1} _{A_1^j} + \infty\mathbf{1}
_{A_2^j}$, and we\vspace*{-3pt} note that $X_{\hat{\tau}}^j = x_1^j$
on the event $\{ \hat{\tau} < \infty\}$.
We are faced with the following possible cases.

If both of $\ell(x_1^j)$, $r(x_1^j)$ are finite and reflecting,
then we appeal to the construction associated with
(\ref{Rlr1})--(\ref{Rlr2}) for $\xi= \xi^j$, $X = X^j$, $\bar{x}
= x_1^j$ and $\tau= \hat{\tau}$\vspace*{-1pt} to obtain processes
$\tilde{\xi}$, $\tilde{X}$ that are equal to $\xi^j$, $X^j$
up to time $\hat{\tau}$ and satisfy (\ref{Rlr2}) for all $t >
\hat{\tau}$.
We then define
\[
\xi^{j+1} = \tilde{\xi} , \qquad X^{j+1} = \tilde{X} \quad
\mbox{and} \quad\tau_{j+1} = \infty\mathbf{1} _{A_1^j} +
\tau_j \mathbf{1} _{A_2^j} .
\]
The result of this construction is such that
$X_{\tau_{j+1}+}^{j+1} = x_2^j$\vspace*{-2pt} on the event
$\{ \tau_{j+1} < \infty\} = A_2^j$, which puts us in the
context of Case~(I) of Step~1.

If $\ell(x_1^j)$ is finite and reflecting and $r(x_1^j) =
\infty$ [resp., $\ell(x_1^j) = -\infty$ and $r(x_1^j)$
is finite and reflecting], then we proceed in the same
way using the construction associated with
(\ref{Rl1})--(\ref{Rl2}) [resp., (\ref{Rr1})--(\ref{Rr2})].

If $\ell(x_1^j)$ is finite and reflecting and $r(x_1^j)$
is finite and repelling, then we consider (\ref{Rl1})--(\ref{Rl2})
and, as above, we construct processes $\tilde{\xi}$,
$\tilde{X}$ that are equal to $\xi^j$, $X^j$ up to time
$\hat{\tau}$ and satisfy (\ref{Rl2}) for all $t > \hat{\tau}$.
We then consider the $(\F_t)$-stopping time $\hat{\tau}
^\dagger$ and the process $\xi^{j+1} \in\acal_\contr$
given by
\begin{eqnarray*}
\hat{\tau} ^\dagger&=& \inf \bigl\{ t \geq\hat{\tau} \mid  \tilde{X}_t \geq r\bigl(x_1^j\bigr) \bigr\}
\quad\mbox{and} \\
 \xi_t^{j+1} &=& %
\cases{
\displaystyle \tilde{\xi}_t ,&\quad if $t \leq\hat{\tau}
^\dagger$,
\cr
\displaystyle \tilde{\xi} _{\hat{\tau} ^\dagger} + \zeta \bigl( r
\bigl(x_1^j\bigr) \bigr) - r\bigl(x_1^j
\bigr) , &\quad if $t > \hat{\tau} ^\dagger$, } %
\end{eqnarray*}
we denote by $X^{j+1}$ the associated solution to
(\ref{X}), and we define
\begin{eqnarray*}
\tau_{j+1} &=& \hat{\tau} ^\dagger\mathbf{1} _{A_1^j} +
\tau_j \mathbf{1} _{A_2^j} ,\qquad    A_1^{j+1}
= \bigl\{ \hat{\tau} ^\dagger< \infty\bigr\} , \qquad  A_2^{j+1}
= A_2^j , \\
   x_i^{j+1} &=&
\zeta \bigl( r\bigl(x_1^j\bigr) \bigr)  \quad  \mbox{and}\quad    x_2^{j+1} = x_2^j
.
\end{eqnarray*}
In this case, we may have $X_{\tau_{j+1}} \in
{\mathcal B}$ but $X_{\tau_{j+1}+} \notin{\mathcal B}$
and $X_t \notin{\mathcal B}$ for all $t \in\,
]\tau_j, \tau_{j+1}[$.

Finally, if $\ell(x_1^j)$ is finite and repelling, and
$r(x_1^j)$ is finite and reflecting, then we are faced
with a construction that is symmetric to the very last
one using (\ref{Rr1})--(\ref{Rr2}).

\emph{STEP 4}.
Go back to Step 2.
\end{pf}

We now prove the main result of the section.
It is worth noting that we can relax significantly assumptions
(\ref{u(X*)-bound})--(\ref{u(x)-bound}).
However, we have opted against any such relaxation because
(a) this would require a considerable amount of extra arguments
of a technical nature that would obscure the main ideas of the
proof, and (b) (\ref{u(X*)-bound})--(\ref{u(x)-bound}) are
plainly satisfied in the special cases that we explicitly solve
in Sections~\ref{12121212} and~\ref{16161616}.

%
\begin{theo} \label{prop:VT}
Consider a function $u\dvtx  \R\rightarrow\R_+$ that satisfies the
conditions of Definition~\ref{u}, let $\xi^* \in\acal_\contr$ be
the control strategy constructed in Lemma~\ref{lem:xi*tau*},
let $X^*$ be the associated solution to (\ref{X}) and define
%
\begin{equation}
v(y) = \max \bigl\{ u(y) , g(y) \bigr\} \quad\mbox{and} \quad w(y) = u(y) \qquad
\mbox{for } y \in\R. \label{v}
\end{equation}
Also, given any $\xi\in\acal_\contr$, define
%
\begin{eqnarray}\label{tau*}
\tau_v^* &=& \tau_v^* (\xi) = \inf \{ t \geq0 \mid X_t \in {\mathcal S} \} ,
\nonumber\\[-8pt]\\[-8pt]\nonumber
\tau_w^* &=&
\tau_w^* (\xi) = \inf \{ t \geq0 \mid X_{t+} \in{\mathcal
S} \} ,
\end{eqnarray}
where $X$ is the associated solution to (\ref{X}), and note
that $\tau_v^* \vee\tau_w^* = \tau_w^*$.
In this context, the following statements are true:
\begin{longlist}[(III)]
\item[(I)]
$J_x^v (\xi^*, \tau) \leq v(x)$ and $J_x^w (\xi^*, \tau)
\leq w(x)$ for all $\tau\in\acal_\stopp$ and all initial
conditions of (\ref{X}).

\item[(II)]
$v(x) = J_x^v (\xi^*, \tau_v^*)$ and $w(x) = J_x^w
(\xi^*, \tau_w^*)$ for every initial condition $x$ of
(\ref{X}) such that
%
\begin{equation}
\sup_{t \geq0} u\bigl(X_t^*\bigr) \leq
K_1   \label{u(X*)-bound}
\end{equation}
for some constant $K_1 = K_1 (x)$.

\item[(III)]
If there exists a constant $K_2$ such that
%
\begin{equation}
u(y) \leq K_2 \qquad\mbox{for all } y \in\R\setminus\scal,
\label{u(x)-bound}
\end{equation}
then $v(x) \leq J_x^v (\xi, \tau_v^*)$ and $w(x) \leq J_x^w
(\xi, \tau_w^*)$ for every initial condition $x$ of~(\ref{X}).

\item[(IV)]
If $u$ satisfies (\ref{u(x)-bound}), then $(\xi^*, \tau_v^*)$
[resp., $(\xi^*, \tau_w^*)$] is an optimal strategy for the game
with performance criterion given by (\ref{Jv}) [resp., (\ref{Jw})]
and $v$ and $w$ are the value functions of the two games.
\end{longlist}
\end{theo}

\begin{pf}
Given a function $u$ satisfying the conditions of
Definition~\ref{u}, we denote by $u''$ the unique,
Lebesgue-a.e., first derivative of $u'$ in $\R\setminus
{\mathcal B}$, and we define $u''(x)$, $u'(x)$ arbitrarily
for $x$ in the finite set ${\mathcal B}$.
In view of (\ref{xi*1}), we can use It\^{o}'s formula and the
integration by parts formula to calculate
\begin{eqnarray*}
e^{- \Lambda_T} u\bigl(X_T^*\bigr) &= & u(x) + \int
_0^T e^{-\Lambda_t} \Lop u
\bigl(X_t^*\bigr) \, dt + \int_{[0, T[}
e^{-\Lambda_t} u'\bigl(X_t^*\bigr) \, d
\xi_t
\\
&&{} + \sum_{0 \leq t < T} e^{- \Lambda_t} \bigl[ u
\bigl(X_{t+}^*\bigr) - u\bigl(X_t^*\bigr) - u'
\bigl(X_t^*\bigr) \Delta X_t^* \bigr] + M_T^*
,
\end{eqnarray*}
where
%
\begin{equation}
M_T^* = \int_0^T e^{-\Lambda_t}
\sigma\bigl(X_t^*\bigr) u'\bigl(X_t^*\bigr)
\, dW_t . \label{M*}
\end{equation}
Rearranging terms and using (\ref{xi*2})--(\ref{xi*3}),
we obtain
\begin{eqnarray*}
&&\int_0^T e^{- \Lambda_t} h
\bigl(X_t^*\bigr) \, dt + \int_{[0, T[}
e^{-\Lambda_t} \, d\check{\xi} _t^* + e^{- \Lambda_T} u
\bigl(X_T^*\bigr)
\\
&&\qquad = u(x) + \int_0^T e^{-\Lambda_t}
\bigl[ \Lop u\bigl(X_t^*\bigr) + h\bigl(X_t^*\bigr) \bigr]
\, dt + \int_0^T e^{-\Lambda_t} \bigl[ 1 +
u'\bigl(X_t^*\bigr) \bigr] \, d \bigl(
\xi^{*{\mathrm c}} \bigr) _t^+
\\
&&\qquad \quad {} + \int_0^T e^{-\Lambda_t}
\bigl[ 1 - u'\bigl(X_t^*\bigr) \bigr] \, d \bigl(
\xi^{*{\mathrm c}} \bigr) _t^- \\
&&\qquad \quad {}+ \sum_{0 \leq t < T}
e^{- \Lambda_t} \bigl[ u\bigl(X_{t+}^*\bigr) - u\bigl(X_t^*
\bigr) + \bigl\llvert \Delta X _t^*\bigr\rrvert \bigr] +
M_T^*
\\
&&\qquad = u(x) + \int_0^T e^{-\Lambda_t}
\bigl[ \Lop u\bigl(X_t^*\bigr) + h\bigl(X_t^*\bigr) \bigr]
\, dt + M_T^* .
\end{eqnarray*}
It follows that, given any finite $(\F_t)$-stopping time
$\hat{\tau}$,
%
\begin{eqnarray}\label{Ito*1}
&&\int_0^{\hat{\tau}} e^{- \Lambda_t} h
\bigl(X_t^*\bigr) \, dt + \int_{[0, \hat{\tau}[}
e^{-\Lambda_t} \, d\check{\xi} _t^* + e^{- \Lambda_{\hat{\tau}}} g
\bigl(X_{\hat{\tau}}^*\bigr)
\nonumber
\\
&&\qquad = u(x) + e^{- \Lambda_{\hat{\tau}}} \bigl[ g\bigl(X_{\hat{\tau}}^*\bigr) - u
\bigl(X_{\hat{\tau}}^*\bigr) \bigr] + \int_0^{\hat{\tau}}
e^{-\Lambda_t} \bigl[ \Lop u\bigl(X_t^*\bigr) + h
\bigl(X_t^*\bigr) \bigr] \, dt + M_{\hat{\tau}}^*
\nonumber
\\[-8pt]
\\[-8pt]
&&\qquad = u(x) \mathbf{1} _{\{ 0 < \hat{\tau} \}} + e^{- \Lambda_{\hat{\tau}}} \bigl[ g
\bigl(X_{\hat{\tau}}^*\bigr) - u\bigl(X_{\hat{\tau}}^*\bigr) \bigr] \mathbf{1}
_{\{ 0 < \hat{\tau} \}} + g(x) \mathbf{1} _{\{ \hat{\tau} = 0 \}}
\nonumber
\\
&&\qquad \quad {} + \int_0^{\hat{\tau}} e^{-\Lambda_t}
\bigl[ \Lop u\bigl(X_t^*\bigr) + h\bigl(X_t^*\bigr)
\bigr] \, dt + M_{\hat{\tau}}^* .
\nonumber
\end{eqnarray}
Similarly, we can calculate
%
\begin{eqnarray}\label{Ito*2}
&&\int_0^{\hat{\tau}} e^{- \Lambda_t} h
\bigl(X_t^*\bigr) \, dt + \int_{[0, \hat{\tau}]}
e^{-\Lambda_t} \, d\check{\xi} _t^* + e^{- \Lambda_{\hat{\tau}}} g
\bigl(X_{\hat{\tau}+}^*\bigr)
\nonumber
\\
&&\qquad = u(x) + e^{- \Lambda_{\hat{\tau}}} \bigl[ g\bigl(X_{\hat{\tau}+}^*\bigr) - u
\bigl(X_{\hat{\tau}+}^*\bigr) \bigr]\\
&&\qquad \quad {} + \int_0^{\hat{\tau}}
e^{-\Lambda_t} \bigl[ \Lop u\bigl(X_t^*\bigr) + h
\bigl(X_t^*\bigr) \bigr] \, dt + M_{\hat{\tau}}^* .
\nonumber
\end{eqnarray}

Combining (\ref{Ito*1}) with (\ref{xi*2}) and
the facts that $0 \leq g(x) \leq u(x)$ for all $x \in\R
\setminus\ccal= \operatorname{int} (\wcal\cup
\scal_\wcal)$ and $\Lop u (x) + h(x) \leq0$
Lebesgue-a.e. in $\R\setminus\ccal$, we can see
that, given any $T>0$ and any $(\F_t)$-stopping
time $\tau$,
%
\begin{eqnarray}
\label{Ito*3}
&&\int_0^{T \wedge\tau} e^{- \Lambda_t} h
\bigl(X_t^*\bigr) \, dt + \int_{[0, T \wedge\tau[}
e^{-\Lambda_t} \, d\check{\xi} _t^* + e^{- \Lambda_\tau} g
\bigl(X_\tau^*\bigr) \mathbf{1 } _{\{ \tau\leq T \}}
\nonumber
\\[-8pt]
\\[-8pt]
&&\qquad \leq u(x) \mathbf{1} _{\{ 0 < \tau\}} + g(x) \mathbf{1}
_{\{ \tau= 0 \}} + M_{T \wedge\tau}^* \leq v(x) + M_{T \wedge\tau}^* ,
\nonumber
\end{eqnarray}
the last inequality following thanks to (\ref{v}).
These inequalities and the positivity of $h$, $g$ imply
that the stopped process
$M^{*\tau}$ is a supermartingale and $\e[M
_{T \wedge\tau}^*] \leq0$.
Therefore, we can take expectations in (\ref{Ito*3})
and pass to the limit $T \rightarrow\infty$ using Fatou's
lemma to obtain the inequality $J_x^v (\xi^*, \tau)
\leq\max \{ u(x) , g(x)  \} = v(x)$.
With reasoning similar to (\ref{Ito*2}), we derive
the inequality $J_x^w (\xi^*, \tau) \leq u(x) = w(x)$,
and (I) follows.

To prove (II), we consider the $(\F_t)$-stopping time
$\tau_v^*$ defined by (\ref{tau*}) with $X^*$ instead of
$X$, and we note that
\[
X_t^* \in\operatorname{cl} \wcal= \R\setminus \operatorname{int} (
\ccal\cup\scal ) \qquad\mbox{for all } 0 < t \leq\tau_v^* .
\]
Combining this observation and the definition of $\tau_v^*$
with the facts that $g(x) \leq u(x) = v(x)$ for all $x \in\wcal$
and $v(x) = g(x)$ for all $x \in\scal$, we can see that
\begin{eqnarray*}
u\bigl(X_{\tau_v^*}^*\bigr) \mathbf{1} _{\{ \tau_v^* > 0 \}} &=& g
\bigl(X_{\tau_v^*}^*\bigr) \mathbf{1} _{\{ \tau_v^* > 0 \}} ,
\\
v(x) \mathbf{1} _{\{ \tau_v^* = 0 \}} &=& g(x) \mathbf{1} _{\{ \tau_v^* = 0 \}} \quad
\mbox{and} \quad v(x) \mathbf{1} _{\{ \tau_v^* > 0 \}} = u(x) \mathbf{1}
_{\{ \tau_v^* > 0 \}} .
\end{eqnarray*}
In view of these observations, (\ref{Ito*1}) and the fact
that $\Lop u(x) + h(x) = 0$ Lebesgue-a.e. in $\wcal$, we can
see that, given any $T>0$,
\begin{eqnarray*}
&&\int_0^{T \wedge\tau_v^*} e^{- \Lambda_t} h
\bigl(X_t^*\bigr) \, dt + \int_{[0, T \wedge\tau_v^*[}
e^{-\Lambda_t} \, d\check{\xi} _t^* + e^{- \Lambda_{\tau_v^*}} g
\bigl(X_{\tau_v^*}^*\bigr) \mathbf{1} _{\{ \tau_v^* \leq T \}}\\
&&  \quad {} + e^{-\Lambda_T} u
\bigl(X_T^*\bigr) \mathbf{1} _{\{ T < \tau_v^* \}}
\\
&&\qquad = u(x) \mathbf{1} _{\{ 0 < \tau_v^* \}} + e^{- \Lambda_{\tau_v^*}} \bigl[ g
\bigl(X_{\tau_v^*}^*\bigr) - u\bigl(X_{\tau_v^*}^*\bigr) \bigr]
\mathbf{1} _{\{ 0 < \tau_v^* \leq T \}} + g(x) \mathbf{1} _{\{ \tau_v^* = 0 \}}\\
&&\qquad \quad {} +
M_{T \wedge\tau_v^*}^*
\\
&&\qquad = v(x) + M_{T \wedge\tau_v^*}^* .
\end{eqnarray*}
If we denote by $(\varrho_n)$ a localizing sequence
for the stopped local martingale $M^{* \tau_v^*}$ such that
$\varrho_n > 0$ for all $n \geq1$, then we can see that
these identities imply that
\begin{eqnarray*}
&&\e \biggl[ \int_0^{\varrho_n \wedge\tau_v^*} e^{- \Lambda_t} h
\bigl(X_t^*\bigr) \, dt + \int_{[0, \varrho_n \wedge\tau_v^*[}
e^{-\Lambda_t} \, d\check{\xi} _t^* + e^{- \Lambda_{\tau_v^*}} g
\bigl(X_{\tau_v^*}^*\bigr) \mathbf{1} _{\{ \tau_v^* \leq\varrho_n \}}
\\
&&\hspace*{185pt}\qquad {} + e^{-\Lambda_{\varrho_n}} u\bigl(X_{\varrho_n}^*\bigr) \mathbf{1}
_{\{ \varrho_n < \tau^* \}} \biggr] = v(x) .
\end{eqnarray*}
In view of (\ref{u(X*)-bound}) and Assumption~\ref{Assm2},
we can pass to the limit as $n \rightarrow\infty$ using the
monotone and the dominated convergence theorems to
obtain $J_x^v (\xi^*, \tau^*) = \max \{ u(x) , g(x)  \}
= v(x)$.

We can use (\ref{Ito*2}) and the observations that
\begin{eqnarray*}
X_t^* \in\operatorname{cl} \wcal= \R\setminus \operatorname{int} (
\ccal\cup\scal )\qquad  \mbox{for all } 0 < t \leq\tau_w^* \quad\mbox{and}
\quad u\bigl(X_{\tau_w^*+}^*\bigr) = g\bigl(X_{\tau_w^*+}^*\bigr)
\end{eqnarray*}
to show that $J_x^w (\xi^*, \tau^*) = u(x) = w(x)$ similarly.

To establish Part~(III), we consider any admissible $\xi\in
\acal_\contr$ and we note that (\ref{Ito*1}) remains true with
$\xi$, $X$ instead of $\xi^*$, $X^*$ if $\hat{\tau}$ is replaced
by $\hat{\tau} \wedge\tau_v^*$ because $\bcal\subseteq
\scal$.
Also, we note that
%
\begin{equation}
X_t \in\R\setminus\scal= ( \wcal\cup\ccal) \setminus \scal\qquad
\mbox{for all } t < \tau_v^* .
\end{equation}
In view of the facts that $g(x) \leq u(x) = v(x)$ for all $x \in
\R\setminus\scal$ and $u(x) \leq g(x) = v(x)$ for all $x \in
\scal$, we can see that this observation and the definition
of $\tau_v^*$ imply that
\begin{eqnarray*}
u(X_{\tau_v^*}) \mathbf{1} _{\{ \tau_v^* > 0 \}} &\leq& g(X_{\tau_v^*})
\mathbf{1} _{\{ \tau_v^* > 0 \}} ,
\\
v(x) \mathbf{1} _{\{ \tau_v^* = 0 \}} &=& g(x) \mathbf{1} _{\{ \tau_v^* = 0 \}} \quad
\mbox{and} \quad v(x) \mathbf{1} _{\{ \tau_v^* > 0 \}} = u(x) \mathbf{1}
_{\{ \tau_v^* > 0 \}} .
\end{eqnarray*}
Combining these observations with the fact that
$\Lop u(x) + h(x) \geq0$ Lebesgue-a.e. inside
$\operatorname{int} [(\wcal\cup\ccal) \setminus\scal]$,
we can see that (\ref{Ito*1}) implies that, given any $T>0$,
\begin{eqnarray*}
&&\int_0^{T \wedge\tau_v^*} e^{- \Lambda_t}
h(X_t) \, dt + \int_{[0, T \wedge\tau_v^*[} e^{-\Lambda_t} \, d
\check{\xi} _t + e^{- \Lambda_{\tau_v^*}} g(X_{\tau_v^*}) \mathbf{1}
_{\{ \tau_v^* \leq T \}} \\
&&  \quad {}+ e^{-\Lambda_T} u(X_T) \mathbf{1}
_{\{ T < \tau_v^* \}}
\\
&&\qquad \geq u(x) \mathbf{1} _{\{ 0 < \tau_v^* \}} + e^{- \Lambda_{\tau_v^*}} \bigl[
g(X_{\tau_v^*}) - u(X_{\tau_v^*}) \bigr] \mathbf{1} _{\{ 0 < \tau_v^* \leq T \}} +
g(x) \mathbf{1} _{\{ \tau_v^* = 0 \}} \\
&&\qquad \quad {}+ M_{T \wedge\tau_v^*}
\\
&&\qquad \geq v(x) + M_{T \wedge\tau_v^*} ,
\end{eqnarray*}
where $M$ is defined as in (\ref{M*}).
If $(\varrho_n)$ is a localizing sequence for the stopped local
martingale $M^{\tau_v^*}$ such that $\varrho_n > 0$ for all
$n \geq1$, then these inequalities imply that
\begin{eqnarray*}
&&\e \biggl[ \int_0^{\varrho_n \wedge\tau_v^*} e^{- \Lambda_t}
h(X_t) \, dt + \int_{[0, \varrho_n \wedge\tau_v^*[} e^{-\Lambda_t} \, d
\check{\xi} _t + e^{- \Lambda_{\tau_v^*}} g(X_{\tau_v^*}) \mathbf{1}
_{\{ \tau_v^* \leq\varrho_n \}}
\\
&&\hspace*{181pt}\qquad{} + e^{-\Lambda_{\varrho_n}} u\bigl(X_{\varrho_n}^*\bigr) \mathbf{1}
_{\{ \varrho_n < \tau^* \}} \biggr] \geq v(x) .
\end{eqnarray*}
In view of (\ref{u(x)-bound}) and Assumption~\ref{Assm2},
we can pass to the limit as $n \rightarrow\infty$ using the
monotone and the dominated convergence theorems to obtain
$J_x^v (\xi, \tau^*) \geq\max \{ u(x) , g(x)  \} = v(x)$.

In general, the inequality $\tau_v^* \leq\tau_w^*$ may be
strict because, for example, we may have $x \in\scal$ and $x +
\Delta\xi_0 \in\R\setminus\scal$.
In such a case, the set $ \{ t \in[0, \tau_w^*[ \, \mid
 X_t \in\bcal \}$ may not be empty, but it is finite.
Therefore, we can use It\^{o}'s formula to derive (\ref{Ito*1})
with $\xi$, $X$ instead of $\xi^*$, $X^*$ and with $\hat{\tau}
\wedge\tau_v^*$ replacing $\hat{\tau}$.
Combining this result with the observations that
\begin{eqnarray*}
X_t \in\operatorname{cl} (\R\setminus\scal) \qquad \mbox{for all } 0 < t <
\tau_w^* \quad\mbox{and} \quad u(X_{\tau_w^*+}) \leq
g(X_{\tau_w^*+}) ,
\end{eqnarray*}
we can derive the inequality $J_x^w (\xi, \tau^*) \geq u(x)
= w(x)$ as above.

Finally, Part~(IV) follows immediately from Parts~(I)--(III).
\end{pf}

%
\begin{rem} \label{rem1}
An inspection of the proof of Theorem~\ref{prop:VT}
reveals that the optimal strategy $(\xi^*, \tau_w^*)$ of the game
where the controller has the first-move advantage is highly
nonunique.
Indeed, in the presence of (\ref{u(x)-bound}), $(\xi^*,
\tilde{\tau}_w^*)$, where $\tilde{\tau}_w^*$ is any
$(\F_t)$-stopping time such that $X_{\tilde{\tau}_w^* +}^*
\mathbf{1} _{\{ \tilde{\tau}_w^* < \infty\}} \in{\mathcal S}$,
in particular, $(\xi^*, \infty)$, is also an optimal strategy.
It is worth noting that a similar observation cannot be made
for the game where the stopper has the first-move
advantage.
Both of the special cases considered in the following two
sections provide cases illustrating this situation; see
Propositions~\ref{prop:caseI3}, \ref{prop:caseI4},
\ref{prop:caseII2} and~\ref{prop:caseII3}.
\end{rem}

\section{The explicit solution to a special case with quadratic
reward functions}\label{sec5}
\label{12121212}

We now derive the explicit solution to the special case of
the general problem that arises when
\begin{eqnarray*}
b(x) &=& 0 , \qquad\sigma(x) = 1, \qquad\bar{\delta} (x) = \delta , \qquad h(x) =
\kappa x^2 + \mu\quad\mbox{and} \\
 g(x)& =& \lambda x^2
\qquad\mbox{for all } x \in\R,
\end{eqnarray*}
for some constants $\delta, \kappa, \lambda> 0$ and
$\mu\geq0$.
In our analysis, we exploit the symmetry around
the origin that the problem has, we consider only
sets $\Gamma\subseteq\R$ such that
$\{ -x \mid x \in\Gamma\} = \Gamma$ and we denote
$\Gamma^+ = \Gamma\cap[0,\infty[$.
Also, we recall that the general solution to the ODE
\[
\Lop f(x) + h(x) \equiv{\tfrac{1}{2}}f''(x) -
\delta f(x) + \kappa x^2 + \mu= 0
\]
is given by
\[
f(x) = A \cosh\sqrt{2\delta} x + B \sinh\sqrt{2\delta} x + \frac{\kappa}{\delta}
x^2 + \frac{\kappa+ \delta\mu} {
\delta^2}
\]
for some constants $A, B \in\R$.

\begin{figure}

\includegraphics{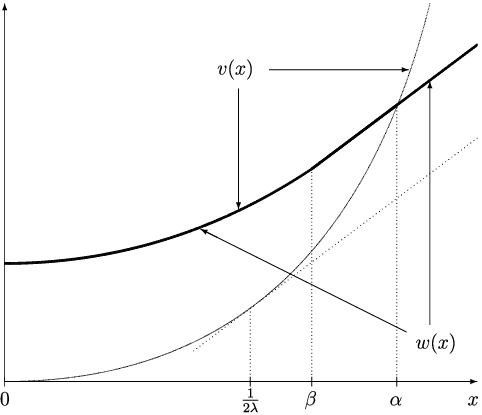}

\caption{The functions $v$
and $w$ in the context of Proposition~\protect\ref{prop:caseI1}
($\beta< \alpha$).}
\label{fig2}\end{figure}

\begin{figure}

\includegraphics{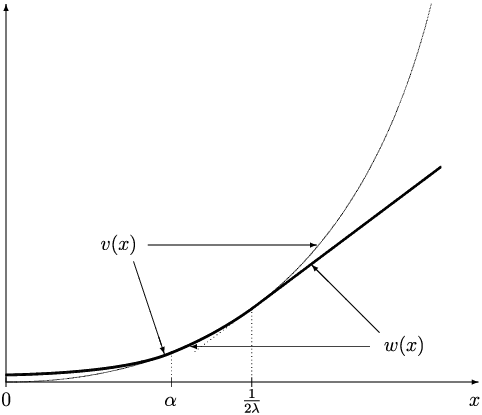}

\caption{The functions $v$
and $w$ in the context of Proposition~\protect\ref{prop:caseI2}
($\alpha< \beta= \frac{1}{2\lambda}$).}
\label{fig3}\end{figure}

In the special case that we consider in this section, the
controller should exert effort to keep the state process
close to the origin.
On the other hand, the stopper should terminate the game
if the state process is sufficiently far from the origin.
In view of these observations, we derive optimal strategies
by considering functions satisfying the requirements of
Definition~\ref{u} that are associated with the regions
%
\begin{equation}
\scal^+ = [\alpha, \infty [ , \qquad \ccal^+ = [\beta, \infty [ \quad\mbox{and}
\quad \wcal^+ = [ 0 , \alpha\wedge\beta [   \label{WSC0}
\end{equation}
for some constants $\alpha, \beta> 0$; see
Definition~\ref{u}.
In particular, we derive three qualitatively different cases
that are characterized by the relations $\beta< \alpha$,
$\alpha< \beta$ or $\alpha= \beta$, depending on
parameter values; see Figures~\ref{fig2}--\ref{fig4} as well as
Remark~\ref{rem:example}.

In this context,
Theorem~\ref{prop:VT} implies that the associated optimal
strategies can be described informally as follows.
The controlled process $\xi^*$ has an initial jump equal
to $-(x+\beta)$ [resp., $-(x-\beta)$] if the initial condition
$x$ of (\ref{X}) is such that $x<-\beta$ (resp., $x>\beta$).
Beyond time 0, $\xi^*$ is such that the associated
solution to (\ref{X}) is reflecting in $-\beta$ in the positive
direction and in $\beta$ in the negative direction.
On the other hand, the optimal stopping times $\tau_v^*$,
$\tau_w^*$ are the first hitting times of ${\mathcal S}$
as defined by (\ref{tau*}).
In view of these observations, we focus on the construction
of the function $u$ satisfying the requirements of
Definition~\ref{u} in what follows.

In the first case that we consider, $u$ identifies with the value
function of the singular stochastic control problem that arises
if the stopper never terminates the game (see Figure~\ref{fig2}).
In particular, we look for a solution to the variational inequality
%
\begin{equation}
\min \bigl\{ {\tfrac{1}{2}}u''(x) - \delta
u(x) + \kappa x^2 + \mu,  1 - \bigl\llvert u'(x)\bigr
\rrvert \bigr\} = 0 \label{u-VI1}
\end{equation}
of the form
%
\begin{equation}
u(x) = %
\cases{\displaystyle A \cosh\sqrt{2\delta} x +
\frac{\kappa}{\delta} x^2 + \frac{\kappa+ \delta\mu}{\delta^2} , &\quad if $|x| \leq
\beta$,
\cr
\displaystyle x - \beta+ u(\beta) , &\quad if $|x| > \beta$. }
\label{u1}
\end{equation}
The requirement that $u$ should be $C^2$ along the
free-boundary point $\beta$, which is associated with the
so-called ``principle of smooth fit'' of singular stochastic
control, implies that the parameter $A$ should be given by
%
\begin{equation}
A = - \frac{\kappa}{\delta^2 \cosh\sqrt{2\delta} \beta} , \label{A1}
\end{equation}
while $\beta> 0$ should satisfy
%
\begin{equation}
\tanh\sqrt{2\delta} \beta= \frac{\delta(2\kappa\beta-
\delta)}{\kappa\sqrt{2\delta}} \label{alpha1} .
\end{equation}
We also define $\alpha> 0$ to be the unique solution to
the equation
%
\begin{equation}
u(\alpha) = \lambda\alpha^2 . \label{alpha-defn1}
\end{equation}


We prove the following result, as well as the other ones
we consider in this section, in Appendix~\ref{app1}.

%
\begin{prop} \label{prop:caseI1}
Equation (\ref{alpha1}) has a unique solution $\beta> 0$,
which is strictly greater than $\frac{\delta}{2\kappa}$,
while equation (\ref{alpha-defn1}) has a unique solution
$\alpha> 0$.
Furthermore, $\alpha> \beta$ if and only if
%
\begin{eqnarray}
\label{cond1}
\quad \delta\lambda- \kappa< 0 \hspace*{5pt}&& \mbox{or} \quad\delta \lambda- \kappa= 0 \quad
\mbox{and}\quad \mu> 0
\nonumber
\\[-4pt]
\\[-12pt]
&&\mbox{or} \quad\delta\lambda- \kappa> 0 \quad \mbox{and}\quad \tanh\sqrt{
\frac{2\delta\mu}{\delta\lambda- \kappa}} < \sqrt{\frac{2\delta\mu}{\delta\lambda- \kappa}} - \frac{\delta^2}{\kappa\sqrt{2\delta}} ,
\nonumber
\end{eqnarray}
in which case, $\alpha> \frac{1}{2\lambda}$ and the
function $u$ defined by (\ref{u1}) for $A<0$, given by
(\ref{A1}), satisfies the conditions of Definition~\ref{u}; see
Figure~\ref{fig2} for a depiction of the value functions
$v$ and $w$.
\end{prop}

%

We next consider the possibility that the value function of the
game where the stopper has the ``first-move advantage'' identifies
with the value function of the optimal stopping problem that
arises if the controller never acts; see Figure~\ref{fig3}.%
In this case, we look for a solution to the variational inequality
\[
\max \bigl\{ {\tfrac{1}{2}}v''(x) - \delta
v(x) + \kappa x^2 + \mu,  \lambda x^2 - v(x) \bigr\} =
0
\]
of the form
%
\begin{equation}
v(x) = %
\cases{\displaystyle A \cosh\sqrt{2\delta} x +
\frac{\kappa}{\delta} x^2 + \frac{\kappa+ \delta\mu}{\delta^2} , &\quad if $|x| \leq \alpha
$,
\cr
\displaystyle \lambda x^2 , &\quad if $|x| > \alpha$. }
\label{w2}
\end{equation}
The requirement that $v$ should be $C^1$ along the
free-boundary point $\alpha$, which is associated with
the so-called ``principle of smooth fit'' of optimal stopping,
implies that the parameter $A$ should be given by
%
\begin{equation}
A = \frac{\delta(\delta\lambda- \kappa) \alpha^2 - (\kappa+
\delta\mu)} {\delta^2 \cosh\sqrt{2\delta} \alpha} , \label{A2}
\end{equation}
while $\alpha> 0$ should satisfy
%
\begin{equation}
\tanh\sqrt{2\delta} \alpha= \frac{\sqrt{2\delta} (\delta\lambda
- \kappa) \alpha}{\delta(\delta\lambda- \kappa) \alpha^2 -
(\kappa+ \delta\mu)} . \label{alpha2}
\end{equation}
In this context, the function $u$ defined by
%
\begin{equation}
u(x) = %
\cases{\displaystyle A \cosh\sqrt{2\delta} x +
\frac{\kappa}{\delta} x^2 + \frac{\kappa+ \delta\mu}{\delta^2} , &\quad if $|x| \leq \alpha
$,
\cr
\displaystyle \lambda x^2 , &\quad if $\displaystyle |x| \in \,\biggl] \alpha,
\frac
{1}{2\lambda} \biggr] $,
\vspace*{2pt}\cr
\displaystyle \frac{1}{4\lambda} , &\quad if $\displaystyle |x| >
\frac
{1}{2\lambda} $, } %
\label{u2}
\end{equation}
provides an appropriate choice for a function satisfying the
requirements of Definition~\ref{u} as long as $\alpha
< \frac{1}{2\lambda}$.



%
\begin{prop} \label{prop:caseI2}
Suppose that $\delta\lambda- \kappa> 0$.
Equation (\ref{alpha2}) has a unique solution $\alpha> 0$,
which is strictly greater than $\sqrt{\frac{\kappa+ \delta\mu} {
\delta(\delta\lambda- \kappa)}}$.\vspace*{-2pt}
This solution is less than or equal to $\frac{1}{2\lambda}$ if
and only if
%
\begin{equation}
\quad \frac{1}{2\lambda} > \sqrt{\frac{\kappa+ \delta\mu} {\delta
(\delta\lambda- \kappa)}} \quad\mbox{and} \quad \tanh
\frac{\sqrt{2\delta}}{2\lambda} \geq\frac{\sqrt{2\delta}
(\delta\lambda- \kappa) \lambda}{\delta(\delta\lambda
- \kappa) - 4 (\kappa+ \delta\mu) \lambda^2} , \label{cond2}
\end{equation}
in which case, the function $u$ defined by (\ref{u2})
for $A>0$, given by (\ref{A2}), satisfies the requirements
of Definition~\ref{u}; see Figure~\ref{fig3} for a depiction of the
value functions $v$ and $w$.
\end{prop}

%

The third case that we consider ``bridges'' the previous
two and is characterized by the fact that the
free-boundary points $\alpha$, $\beta$ may coincide
in a generic way.
In particular, we look for a function $u$ satisfying the
requirements of Definition~\ref{u} that is given by
%
\begin{equation}
u(x) = %
\cases{\displaystyle A \cosh\sqrt{2\delta} x +
\frac{\kappa}{\delta} x^2 + \frac{\kappa+ \delta\mu}{\delta^2} , &\quad if $|x| \leq \alpha
$,
\cr
\displaystyle x - \alpha+ u(\alpha), &\quad if $|x| > \alpha$, }
\label{u3}
\end{equation}
for some $\alpha>0$, and satisfies
%
\begin{equation}
u(\alpha) = \lambda\alpha^2 \label{alpha-defn2};
\end{equation}
see Figure~\ref{fig4}.
The requirements that $u$ should satisfy (\ref{alpha-defn2})
and be $C^1$ at $\alpha$ imply that the parameter $A$
should be given by
%
\begin{equation}
A = \frac{\delta(\delta\lambda- \kappa) \alpha^2 - (\kappa
+ \delta\mu)} {\delta^2 \cosh\sqrt{2\delta} \alpha} , \label{A3}
\end{equation}
while the free-boundary point $\alpha> 0$ should satisfy
%
\begin{equation}
\tanh\sqrt{2\delta} \alpha= \frac{\delta(\delta- 2\kappa
\alpha)} {\sqrt{2\delta}  [ \delta(\delta\lambda- \kappa)
\alpha^2 - (\kappa+ \delta\mu)  ]} . \label{alpha3}
\end{equation}


%
\begin{prop} \label{prop:caseI3}
Suppose that $\delta\lambda- \kappa> 0$ and $\sqrt{
\frac{\kappa+ \delta\mu} {\delta(\delta\lambda- \kappa)}}
\neq\frac{\delta}{2\kappa}$.%
\begin{figure}

\includegraphics{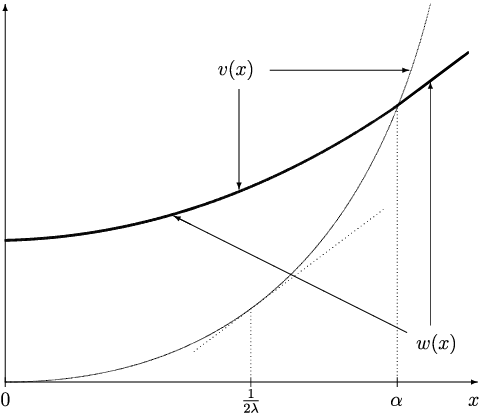}

\caption{The functions $v$
and $w$ in the context of Propositions~\protect\ref{prop:caseI3}
and~\protect\ref{prop:caseI4} ($\alpha= \beta$).}
\label{fig4}\end{figure}
Equation (\ref{alpha3}) has a unique solution $\alpha> 0$
such that
%
\begin{eqnarray}
\mbox{if} \quad\frac{\delta}{2\kappa} < \sqrt{\frac{\kappa+
\delta
\mu} {\delta(\delta\lambda- \kappa)}} , \quad\mbox{then}
\quad \frac{1}{2\lambda} < \frac{\delta}{2\kappa} < \alpha< \sqrt{\frac{\kappa+ \delta\mu} {\delta(\delta\lambda- \kappa)}} ,
\label{alpha3range1}
\end{eqnarray}
while
\begin{eqnarray}
\mbox{if} \quad\sqrt{\frac{\kappa+ \delta\mu} {\delta(\delta
\lambda- \kappa)}} < \frac{\delta}{2\kappa} , \quad\mbox{then}
\quad\sqrt{\frac{\kappa+ \delta\mu} {\delta(\delta\lambda-
\kappa)}} < \alpha< \frac{\delta}{2\kappa} . \label{alpha3range2}
\end{eqnarray}
If the parameters are such that (\ref{alpha3range1}) is true,
then the function $u$ defined by (\ref{u3}) for $A<0$, given
by (\ref{A3}), satisfies the conditions of Definition~\ref{u}
if and only if
%
\begin{equation}
\tanh\sqrt{\frac{2\delta\mu}{\delta\lambda- \kappa}} \geq \sqrt{\frac{2\delta\mu}{\delta\lambda- \kappa}} -
\frac{\delta^2}{\kappa\sqrt{2\delta}} . \label{cond31}
\end{equation}
On the other hand, if the parameters are such that
(\ref{alpha3range2}) is true, then $\frac{1}{2\lambda} < \alpha$
if and only if
%
\begin{eqnarray}\label{cond32}
\frac{1}{2\lambda} & \leq&\sqrt{\frac{\kappa+ \delta\mu} {
\delta(\delta\lambda- \kappa)}} \quad \mbox{or}
\\
\frac{1}{2\lambda} &>& \sqrt{\frac{\kappa+
\delta\mu} {\delta(\delta\lambda- \kappa)}} \quad\mbox{and} \quad\tanh
\frac{\sqrt{2\delta}}{2\lambda} < \frac{\sqrt{2\delta}
(\delta\lambda- \kappa) \lambda}{\delta(\delta\lambda
- \kappa) - 4 (\kappa+ \delta\mu) \lambda^2 } ,
\nonumber
\end{eqnarray}
in which case, the function $u$ defined by (\ref{u3}) for
$A>0$, given by (\ref{A3}), satisfies the conditions of
Definition~\ref{u}; see Figure~\ref{fig4} for a depiction of the
value functions $v$ and~$w$.
\end{prop}

%

The results that we have established thus far involve
mutually exclusive conditions on the problem data.
To exhaust all possible parameter values, we need to
consider the following result that is associated with
the regions
%
\begin{equation}
\bcal= \scal_\wcal= \varnothing, \qquad\ccal^+ = \scal_\ccal^+
= \biggl[ \frac{\delta}{2\kappa}, \infty \biggr[ \quad\mbox{and} \quad\wcal^+ = \biggl[ 0,
\frac{\delta}{2\kappa} \biggr[ ,
\end{equation}
which are consistent with (\ref{WSC0}) for $\alpha=
\beta= \frac{\delta}{2\kappa}$, and the proof of which
is straightforward.

%
\begin{prop} \label{prop:caseI4}
Suppose that $\delta\lambda- \kappa> 0$ and $\sqrt{
\frac{\kappa+ \delta\mu} {\delta(\delta\lambda- \kappa)}}
= \frac{\delta}{2\kappa}$.
The function $u$ defined by
%
\begin{equation}
u(x) = %
\cases{\displaystyle \frac{\kappa}{\delta} x^2 +
\frac{\kappa
+ \delta\mu}{\delta^2} , &\quad if $\displaystyle |x| \leq\frac{\delta} {
2\kappa} $,
\cr
\displaystyle x
- \frac{\delta}{2\kappa} + \frac{\lambda
\delta^2}{4\kappa^2} , &\quad if $\displaystyle |x| > \frac{\delta} {2\kappa} $,
} %
\end{equation}
is a $C^1$ function that satisfies the requirements of
Definition~\ref{u}.
\end{prop}

%
\begin{rem} \label{rem:example}
Suppose that $\delta\lambda- \kappa> 0$.
The conditions differentiating between the different cases
we have considered are mutually exclusive and exhaustive
in the sense that they cover the entire range of possible
parameter values.
To see this claim, we define
%
\[
Q_1 = \sqrt{\frac{2\delta\mu}{\delta\lambda- \kappa}} - \frac{\delta^2}{\kappa\sqrt{2\delta}} , \qquad
Q_2 = \sqrt{\frac{\kappa+ \delta\mu} {\delta(\delta\lambda
- \kappa)}}
\]
{and}
\begin{equation}
Q_3 = \frac{\sqrt{2\delta} (\delta\lambda- \kappa) \lambda} {
\delta(\delta\lambda- \kappa) - 4 (\kappa+ \delta\mu)
\lambda^2} =\frac{2\lambda}{\sqrt{2\delta}  ( 1 - 4
\lambda^2 Q_2^2  )} .
\end{equation}
In view of the implications
\[
Q_1 > 0  \quad \Rightarrow\quad  \frac{\delta}{2\kappa} < Q_2 \quad
\mbox{and} \quad Q_2 < \frac{1}{2\lambda}  \quad \Leftrightarrow\quad  0 <
Q_3 ,
\]
we can see that the following table summarizes the conditions
of Propositions~\ref{prop:caseI1}, \ref{prop:caseI2},
\ref{prop:caseI3} and~\ref{prop:caseI4}:\vspace*{6pt}
\begin{center}
\begin{tabular}{@{}ll@{}}
\hline
Proposition~\ref{prop:caseI1} $ ( \beta< \alpha )$
& $\tanh\sqrt{\frac{2\delta\mu}{\delta\lambda- \kappa}}
< Q_1$ \\
Proposition~\ref{prop:caseI2} $ ( \alpha< \beta=
\frac{1}{2\lambda}  )$
& $\frac{1}{2\lambda} > Q_2$ and $\tanh
\frac{\sqrt{2\delta}}{2\lambda} > Q_3$ \\[3pt]
Propositions~\ref{prop:caseI3}, \ref{prop:caseI4}
($\beta= \alpha$) & $\frac{\delta}{2\kappa} < Q_2$ and
$\tanh\sqrt{\frac{2\delta\mu}{\delta\lambda- \kappa}} > Q_1$ \\
& \emph{or} $\frac{1}{2\lambda} \leq Q_2$ \\
& \emph{or} $\frac{1}{2\lambda} > Q_2$ and $\tanh
\frac{\sqrt{2\delta}}{2\lambda} < Q_3$
\\
\hline
\end{tabular}
\end{center}\vspace*{0.5\baselineskip}
For instance, if
\[
\delta= 4 , \qquad\kappa= 1 , \qquad\lambda= \tfrac{1}{2} \quad\mbox{and}
\quad\mu= 9 ,
\]
then $Q_1 = 2 \sqrt{2}$, and we are in the context of
Proposition~\ref{prop:caseI1} if
\[
\delta= 2 , \qquad\kappa= \tfrac{1}{100} , \qquad \lambda=
\tfrac{1}{2} \quad\mbox{and} \quad\mu= 0 ,
\]
then $\frac{1}{2\lambda} = 1 > \frac{1}{\sqrt{198}} = Q_2$,
$\tanh\frac{\sqrt{2\delta}}{2\lambda} \simeq0.9640
> 0.5025 \simeq\frac{99}{197} = Q_3$, and we are in the
context of Proposition~\ref{prop:caseI2} if
\[
\delta= 2 , \qquad\kappa= \tfrac{199}{300} , \qquad \lambda=
\tfrac{1}{2} \quad\mbox{and} \quad\mu= 0 ,
\]
then $\frac{1}{2\lambda} = 1 > \sqrt{\frac{199}{202}} = Q_2$,
$\tanh\frac{\sqrt{2\delta}}{2\lambda} = \tanh2 <
\frac{101}{3} = Q_3$, and we are in the context
of Proposition~\ref{prop:caseI3}, while if
\[
\delta= \tfrac{1}{2} , \qquad\kappa= \tfrac{1}{8} , \qquad \lambda=
\tfrac{1}{2} \quad\mbox{and} \quad\mu= 0 ,
\]
then $\frac{1}{2\lambda} = 1 < \sqrt{2} = Q_2$, and we
are again in the context of Proposition~\ref{prop:caseI3}.
\end{rem}

\section{A special case with value functions that are not \texorpdfstring{$\bolds{C^1}$}{$C^1$}}\label{sec6}
\label{16161616}

We now solve the special case of the general problem that arises
when
\begin{eqnarray*}
b &\equiv&0 , \qquad\sigma\equiv1, \qquad\bar{\delta} \equiv \delta, \qquad h
\equiv0 \quad\mbox{and} \\
 g(x) &=& %
\cases{\displaystyle - \lambda
x^2 + \lambda, &\quad if $|x| \in[0,1] $,
\cr
\displaystyle 0 , &\quad
if $|x| > 1 $, } %
\end{eqnarray*}
for some constants $\delta, \lambda> 0$.
In this context, the controller has no incentive to exert any
control action other than to counter the stopper's action
because \mbox{$h \equiv0$}.
We therefore solve the problem by first viewing the game
from the stopper's perspective.
Also, we exploit the problem's symmetry around the origin
in the same way as in the previous section.

\begin{figure}

\includegraphics{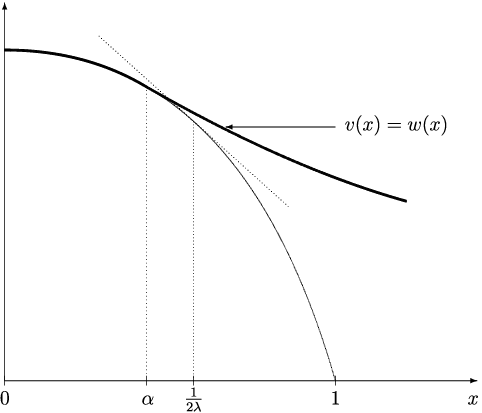}

\caption{The functions $v$
and $w$ in the context of Proposition~\protect\ref{prop:caseII1}.}
\label{fig5}\end{figure}

We first consider the possibility that a function $u$ satisfying
the requirements of Definition~\ref{u} identifies with the value
function of the optimal stopping problem that arises if the
controller never takes any action.
To this end, we look for a solution to the variational inequality
\[
\max \bigl\{ {\tfrac{1}{2}}u''(x) - \delta
u(x) ,  - \lambda x^2 + \lambda - u(x) \bigr\} = 0
\]
of the form
%
\begin{equation}
u(x) = %
\cases{\displaystyle - \lambda x^2 + \lambda, &
\quad if $|x| \leq\alpha$,
\cr
\displaystyle A e^{-\sqrt{2\delta} x} , &\quad if $|x| >
\alpha $, } %
\label{u21}
\end{equation}
for some constants $A$ and $\alpha\in\, ]0,1[$.
A function of this form is associated with the regions
%
\begin{equation}
\bcal= \ccal= \scal_\ccal= \varnothing, \qquad\scal_\wcal^+ =
[0, \alpha] \quad\mbox{and} \quad\wcal^+ = \, ]\alpha, \infty[ ,
\label{WSC21}
\end{equation}
and is depicted by Figure~\ref{fig5}.
To determine the constant $A$ and the free-boundary point
$\alpha$, we appeal to the so-called ``principle of smooth-fit''
of optimal stopping.
We therefore require that $u$ is $C^1$ at $-\alpha$ and
$\alpha$ to obtain
%
\begin{equation}
A = \lambda\bigl(1 - \alpha^2\bigr) e^{\sqrt{2\delta} \alpha} \quad\mbox{and}
\quad \alpha= - \frac{1}{\sqrt{2\delta}} + \sqrt{\frac{1}{2\delta} +1} .
\label{Aalpha21}
\end{equation}

In this case, Theorem~\ref{prop:VT} implies that the
associated optimal strategy can be described informally
as follows.
The controller should never act (i.e., $\xi^* = 0$), while the
stopper should terminate the game as soon as the state
process takes values in ${\mathcal S} = [-\alpha, \alpha]$
(i.e., $\tau_v^* = \tau_u^*$ is the first hitting time of
$[-\alpha, \alpha]$).

%

We prove the following result, as well as the other ones
we consider in this section, in Appendix~\ref{app2}.

%
\begin{prop} \label{prop:caseII1}
The function $u$ defined by (\ref{u21}) for $A > 0$, $\alpha
\in\, ]0,1[$ given by (\ref{Aalpha21}) satisfies the
requirements of Definition~\ref{u} if and only if
%
\begin{equation}
\alpha\leq\frac{1}{2\lambda} \quad\Leftrightarrow\quad \lambda\leq
\frac{1}{2} \biggl( - \frac{1}{\sqrt{2\delta}} + \sqrt{\frac{1}{2\delta} +1}
\biggr) ^{-1}; \label{cond21}
\end{equation}
see Figure~\ref{fig5} for a depiction of the value functions
$v$ and $w$.
\end{prop}

\begin{figure}

\includegraphics{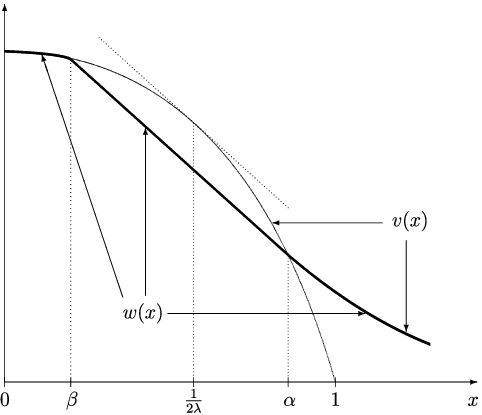}

\caption{The functions $v$
and $w$ in the context of Proposition~\protect\ref{prop:caseII2}.}
\label{fig6}\end{figure}

If the problem data is such that (\ref{cond21}) is not true,
then we consider the possibility that an optimal
strategy is characterized by a function $u$ satisfying the
requirements of Definition~\ref{u} that is associated
with the regions
%
\begin{eqnarray}
\bcal^+ &=& \{ \beta\} , \qquad\scal_\wcal^+ = [0, \beta] , \qquad
\ccal^+ = \scal_\ccal^+ = [\beta, \alpha] \quad \mbox{and} \nonumber
\\[-8pt]
\\[-8pt]\wcal^+
&=& \, ]\alpha, \infty[   \label{WSC22}
\nonumber
\end{eqnarray}
for some $0 \leq\beta< \alpha< 1$, and is depicted by
Figure~\ref{fig6}.
In particular, we consider the function
%
\begin{equation}
u(x) = %
\cases{\displaystyle - \lambda x^2 + \lambda, &
\quad if $|x| \leq\beta$,
\cr
\displaystyle - x - \lambda\alpha^2 +
\alpha+ \lambda , &\quad if $|x| \in \,]\beta, \alpha] $,
\cr
\displaystyle A
e^{-\sqrt
{2\delta} x} , &\quad if $|x| > \alpha$. } %
\label{u22}
\end{equation}
The requirement that $u$ should be continuous at $\beta$
yields
%
\begin{equation}
\lambda\beta^2 - \beta= \lambda\alpha^2 - \alpha,
\label{beta22}
\end{equation}
while, the requirement that $u$ should be $C^1$ along
$-\alpha$, $\alpha$, implies that
%
\begin{equation}
A = \lambda\bigl(1 - \alpha^2\bigr) e^{\sqrt{2\delta} \alpha} \quad\mbox{and}
\quad \alpha= \sqrt{1 - \frac{1}{\lambda\sqrt{2\delta}}} . \label{Aalpha22}
\end{equation}

%

In view of Theorem~\ref{prop:VT}, we can describe
informally the associated optimal strategy as follows.
If the initial condition $x$ of (\ref{X}) belongs to
$ ] {-}\beta, \beta [$, then the controller should wait
until the uncontrolled state process hits $\{ -\beta, \beta
\}$, at which time, the controller should apply an impulse
to instantaneously reposition the state process at
$-\alpha$ or $\alpha$, whichever point is closest.
As soon as the state process takes values in $ ]
{-}\infty, -\alpha ]$ (resp., $[\alpha, \infty[$), the
controller should exert minimal effort to reflect the state
process in $-\alpha$ in the negative direction (resp., in
$\alpha$ in the positive direction).%
\begin{figure}[b]

\includegraphics{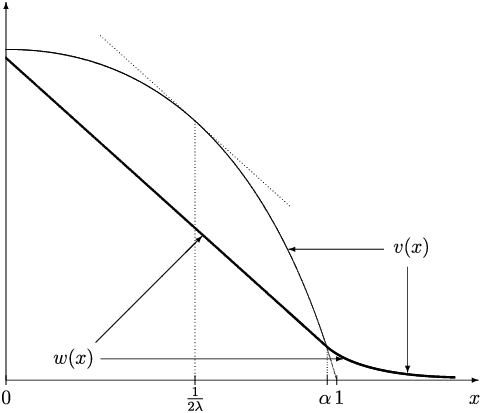}

\caption{The functions $v$
and $w$ in the context of Proposition~\protect\ref{prop:caseII3}.}
\label{fig7}\end{figure}
On the other hand, the stopper should terminate
the game as soon as the state process takes values
in ${\mathcal S} = [-\alpha, \alpha]$.

%
\begin{prop} \label{prop:caseII2}
The point $\alpha$ defined by (\ref{Aalpha22}) is
strictly greater than $\frac{1}{2\lambda}$, and there
exists $\beta\in[0,\alpha[$ satisfying (\ref{beta22})
if and only if
%
\begin{equation}
\frac{1}{2} \biggl( - \frac{1}{\sqrt{2\delta}} + \sqrt{\frac
{1}{2\delta} + 1}
\biggr) ^{-1} < \lambda\leq \biggl( - \frac{1}{\sqrt{2\delta
}} + \sqrt{
\frac{1}{8\delta} + 1} \biggr) ^{-1} , \label{cond22}
\end{equation}
in which case, $\beta< \frac{1}{2\lambda}$.
If the problem data satisfy these inequalities, then the
function $u$ defined by (\ref{u22}), for $A > 0$, $\alpha\in
\, ]0,1[$ given by (\ref{Aalpha22}), satisfies the
conditions of Definition~\ref{u}; see Figure~\ref{fig6} for a
depiction of the value functions $v$ and $w$.
\end{prop}

The final possibility that may arise is associated with the
regions
%
\begin{equation}
\quad \bcal= \{ 0 \} , \qquad\scal_\wcal= \varnothing, \qquad\ccal^+ =
\scal_\ccal^+ = [0, \alpha] \quad \mbox{and} \quad\wcal^+ = \, ]
\alpha, \infty[   \label{WSC23}
\end{equation}
for some $\alpha\in\, ]0,1[$, and is depicted by
Figure~\ref{fig7}.
In this case, a function $u$ satisfying the requirements of
Definition~\ref{u} is given by
%
\begin{equation}
u(x) = %
\cases{\displaystyle - x - \lambda\alpha^2 +
\alpha+ \lambda, &\quad if $|x| \in[0, \alpha] $,
\cr
\displaystyle A
e^{-\sqrt{2\delta
} x} ,&\quad if $|x| > \alpha$. } %
\label{u23}
\end{equation}
The constant $A$ and the free-boundary point $\alpha$
are characterized by the requirement that $u$ should be $C^1$
along $-\alpha$, $\alpha$, and are given by (\ref{Aalpha22}).

In this case, Theorem~\ref{prop:VT} implies that the
associated optimal strategy can be described informally
as follows.
The controlled process $\xi^*$ has an initial jump equal
to $-(x+\alpha)$ [resp., $-(x-\alpha)$] if the initial condition
$x$ of (\ref{X}) is such that $x \in \,] {-}\alpha, 0  ]$
(resp., $x \in\, ]0, \alpha]$).
Beyond time 0, $\xi^*$ is such that the associated
solution to (\ref{X}) is reflecting in $-\alpha$ in the negative
direction if $X_{0+}^* \leq-\alpha$ and in $\alpha$ in the
positive direction if $X_{0+}^* \geq\alpha$.
On the other hand, the stopping time $\tau_v^* =
\tau_v^*$ is the first hitting time of ${\mathcal S}
= [-\alpha, \alpha]$.

%
\begin{prop} \label{prop:caseII3}
The function $u$ defined by (\ref{u22}) for $A > 0$,
$\alpha\in\, ]0,1[$ given by (\ref{Aalpha22})
satisfies the conditions of Definition~\ref{u}
if and only if
%
\begin{equation}
\biggl( - \frac{1}{\sqrt{2\delta}} + \sqrt{\frac{1}{8\delta} + 1} \biggr)
^{-1} < \lambda\label{cond23};
\end{equation}
see Figure~\ref{fig7} for a depiction of the value functions
$v$ and $w$.
\end{prop}

\begin{appendix}
\renewcommand{\thesection}{\Roman{section}}
\section{Proofs of results in Section~\texorpdfstring{\protect\ref{12121212}}{5}}\label{app1}

\begin{pf*}{Proof of Proposition~\texorpdfstring{\ref{prop:caseI1}}{2}}
It is straightforward to see that equation (\ref{alpha1}) has
a unique solution $\beta> 0$ and that this solution is strictly
greater than $\frac{\delta}{2\kappa}$.
In particular, we can verify that
%
\begin{equation}
\tanh\sqrt{2\delta} x - \frac{\delta(2\kappa x - \delta)} {
\kappa\sqrt{2\delta}} %
\cases{\displaystyle > 0,
&\quad for all $x \in [0,\beta[ $,
\cr
\displaystyle < 0 , & \quad for all $x \in\,
]\beta, \infty[ $. } %
\label{tanh-alpha-ineqs1}
\end{equation}
For this value of $\beta$ and for $A<0$ given by (\ref{A1}),
the function $u$ defined by (\ref{u1}) is $C^2$ and satisfies
the variational inequality (\ref{u-VI1}) because
%
\begin{eqnarray}
\bigl\llvert u'(x)\bigr\rrvert &\leq&1 \qquad\mbox{for all }
\llvert x\rrvert \in[0, \beta] , \label{uVI1}
\\
\label{uVI2} \Lop u(x) + h(x) &\equiv&{\tfrac{1}{2}}u''(x)
- \delta u(x) + \kappa x^2 + \mu
\nonumber
\\[-8pt]
\\[-8pt]
&\geq&0 \qquad \mbox{for all } \llvert x\rrvert \in[\beta, \infty[ .
\nonumber
\end{eqnarray}
To see (\ref{uVI1}), we first note that $u'''(x) = (2\delta)
^{\fraca{3}{2}} A \sinh\sqrt{2\delta} x < 0$ for all
$x \in[0, \beta[$, which implies that the restriction of $u''$
in $[0, \beta]$ is strictly decreasing.
Combining this observation with the identities
\[
u''(0) = \frac{2\kappa}{\delta} \biggl( 1 -
\frac{1}{\cosh\sqrt {2\delta}
\beta} \biggr) > 0 \quad\mbox{and} \quad u''(
\beta) = 0 ,
\]
we can see that $u''(x) > 0$ for all $x \in[0, \beta[$.
It follows that $u$ is an even convex function, which, combined
with the identities $u'(0) = 0$ and $u'(\beta) = 1$, implies~(\ref{uVI1}).

To prove (\ref{uVI2}), it suffices to show that
%
\begin{equation}
f_0 (x) \geq0 \qquad\mbox{for all } x \geq\beta, \label{uVI3}
\end{equation}
where
\[
f_0 (x) = {\tfrac{1}{2}}u''(x) -
\delta u(x) + \kappa x^2 + \mu= \kappa x^2 - \delta x +
\delta\beta- \delta u(\beta) + \mu.
\]
The definition and the $C^2$ continuity of $u$ imply that
$f_0 (\beta) = 0$, for $x \geq0$.
Combining this observation with the inequality $f_0' (x) = 2\kappa
( x - \frac{\delta}{2\kappa}  ) > 0$\vspace*{-1pt} for all $x \geq
\beta$,
which follows from the fact that $\beta> \frac{\delta}{2\kappa}$,
we can see that (\ref{uVI3}) is true.\looseness=-1\vadjust{\goodbreak}

To see that equation (\ref{alpha-defn1}) has a unique solution
$\alpha> 0$, we define $f_1 (x) = \lambda x^2 - u(x)$.
In view of the calculations
\begin{eqnarray*}
f_1''' (x) &=& - (2\delta)
^{\fraca{3}{2}} A \sinh\sqrt{2\delta} x > 0  \qquad  \mbox{for } x < \beta
\quad\mbox{and} \\
 f_1'' (x) &=& 2\lambda
x > 0 \qquad  \mbox{for } x > \beta,
\end{eqnarray*}
we can see that either $f_1$ is convex, or there exists
$x_1 \in\, ]0, \beta[$ such that $f_1''(x) < 0$ for
all $x < x_1$ and $f_1''(x) > 0$ for all $x > x_1$.
In the first case, $f_1' (x) > 0$ for all $x > 0$, while,
in the second case, there exists $x_2 > x_1$ such that
$f_1' (x) < 0$ for all $x \in\, ]0, x_2[$ and
$f_1' (x) > 0$ for all $x > x_2$ because $f_1' (0) = 0$.
In either case, we can see that the equation $f_1(x)
= 0$ has a unique solution $\alpha> 0$ because
\[
f_1 (0) = - \frac{\kappa}{\delta^2} \biggl( 1 - \frac{1}{\cosh\sqrt{2\delta} \beta}
\biggr) - \frac{\mu}{\delta} < 0 \quad\mbox{and} \quad \lim_{x \rightarrow\infty}
f_1 (x) = \infty.
\]

To show that the point $\alpha$ defined by (\ref{alpha-defn1})
is strictly greater than $\beta$ if and only if (\ref{cond1}) is true,
we note that the linearity of $u$ in $[\beta, \infty[$ implies that
there exists $\alpha> \beta$ such that (\ref{alpha-defn1}) is
true if and only if $u(\beta) > \lambda\beta^2$.
In particular, if such $\alpha$ exists, then $\alpha> \frac{1} {
2\lambda}$.
Using the definition (\ref{u1}) of $u$, we calculate
\[
u(x) - \lambda x^2 = \frac{\kappa}{\delta^2} \biggl( 1 -
\frac{\cosh\sqrt{2\delta} x}{\cosh\sqrt{2\delta} \beta} \biggr) - \frac{\delta\lambda- \kappa}{\delta} x^2 +
\frac{\mu}{\delta} \qquad\mbox{for } \llvert x\rrvert \leq\beta.
\]
If $\delta\lambda- \kappa< 0$, then this identity implies trivially
that
%
\begin{equation}
u(x) > \lambda x^2 \qquad\mbox{for all } \llvert x\rrvert \leq\beta.
\label{uVI4}
\end{equation}
Similarly, if $\delta\lambda- \kappa= 0$ and $\mu> 0$, then
(\ref{uVI4}) is true.
On the other hand, if $\delta\lambda- \kappa> 0$,
then (\ref{uVI4}) is true if and only if $\beta< \sqrt{\frac{\mu} {
\delta\lambda- \kappa}}$ because the function $x \mapsto
u(x) - \lambda x^2$ is strictly decreasing in $[0,\beta]$.
Therefore, if $\delta\lambda- \kappa\geq0$, then (\ref{uVI4})
is true if and only if the very last inequality in (\ref{cond1}) holds
true, thanks to (\ref{tanh-alpha-ineqs1}).
It follows that the equation $u(x) = \lambda x^2$ has a unique
solution $\alpha> \beta\vee\frac{1}{2\lambda}$ if and only if
(\ref{cond1}) is true.

Finally, it is straightforward to check that, if (\ref{cond1}) is
true, then $u$ is associated with the regions $\bcal=
\scal_\wcal= \varnothing$, $\ccal^+ = [\beta, \infty[$,
$\scal_\ccal^+ = [\alpha, \infty[$ and $\wcal^+ = [0, \beta[$,
and satisfies all of the conditions required by
Definition~\ref{u}.
\end{pf*}

\begin{pf*}{Proof of Proposition~\ref{prop:caseI2}}
The calculation
\[
\frac{d}{d\alpha} \frac{\alpha}{\delta(\delta\lambda- \kappa)
\alpha^2 - (\kappa+ \delta\mu)} = - \frac{\delta(\delta\lambda
- \kappa) \alpha^2 + \kappa+ \delta\mu}{ [ \delta(\delta
\lambda- \kappa) \alpha^2 - (\kappa+ \delta\mu)  ] ^2} < 0
\]
implies that the right-hand side of (\ref{alpha2}) defines a strictly\vspace*{-1pt}
decreasing function on $\R_+ \setminus \{ \sqrt{\frac{\kappa
+ \delta\mu}{\delta(\delta\lambda- \kappa)}}  \}$.
Combining this observation with the fact that $\tanh$ is a strictly
increasing function, we can see that (\ref{alpha2}) has a
unique solution $\alpha> 0$ and that this solution is strictly
greater than $\sqrt{\frac{\kappa+ \delta\mu} {\delta(\delta
\lambda- \kappa)}}$.
In particular, we can see that
%
\begin{eqnarray}
&&\tanh\sqrt{2\delta} x - \frac{\sqrt{2\delta} (\delta\lambda
- \kappa) x}{\delta(\delta\lambda- \kappa) x^2 - (\kappa
+ \delta\mu)} \nonumber
\\[-8pt]
\\[-8pt]
&&\qquad \cases{ > 0 , &\quad if
$\displaystyle x \in \,\biggl] 0, \sqrt{\frac{\kappa+ \delta\mu} {\delta(\delta\lambda-
\kappa)}} \biggr[ \,\cup\,  ]\alpha, \infty[ $,
\vspace*{2pt}\cr
< 0 ,
&\quad if $\displaystyle x \in \,\biggl] \sqrt{\frac{\kappa+ \delta\mu} {
\delta(\delta\lambda- \kappa)}} , \alpha\biggr [ $, } %
\label{tanh-alpha-ineqs2}
\nonumber
\end{eqnarray}
which implies that the solution $\alpha$ of (\ref{alpha2}) is
less than or equal to $\frac{1}{2\lambda}$ if and only if the
inequalities in (\ref{cond2}) are true.

In what follows, we assume that the problem data satisfy
(\ref{cond2}), in which case, $u$ is associated with the
regions $\bcal= \varnothing$, $\scal_\wcal^+ =  [ \alpha,
\frac{1}{2\lambda}  ]$, $\ccal^+ = \scal_\ccal^+ =
[ \frac{1}{2\lambda}, \infty [$ and $\wcal^+ =
[0, \alpha[$.
We will show that $u$ satisfies all of the conditions in
Definition~\ref{u} if and only if we prove that
%
\begin{eqnarray}
u'(x) &\leq&1 \qquad \mbox{for all } x \in[0, \alpha] ,
\label{w2VI1}
\\
u(x) - \lambda x^2 &\geq&0 \qquad \mbox{for all } x \in[0, \alpha] ,
\label{w2VI2}
\\
\label{w2VI3} 
\Lop u(x) + h(x) &\equiv& \frac{1}{2}
u''(x) - \delta u(x) + \kappa x^2 + \mu
\nonumber
\\[-8pt]
\\[-8pt]
&\leq&0 \qquad \mbox{for all } x \in \,\biggl] \alpha, \frac{1}{2\lambda} \biggr[ .
\nonumber
\end{eqnarray}
Inequality (\ref{w2VI1}) follows immediately from the
convexity of $u$ and the fact that $u'(\alpha) = 2\lambda
\alpha\leq1$.
Inequality (\ref{w2VI2}) is equivalent to
%
\begin{equation}
\frac{\delta(\delta\lambda- \kappa) \alpha^2 - (\kappa+
\delta\mu)} {\cosh\sqrt{2\delta} \alpha} \geq f_2 (x) \qquad \mbox{for all } x \in[0,
\alpha] , \label{w2VI4}
\end{equation}
where
\[
f_2 (x) = \frac{\delta(\delta\lambda- \kappa) x^2 -
(\kappa+ \delta\mu)} {\cosh\sqrt{2\delta} x}.
\]
Since $\alpha> \sqrt{\frac{\kappa+ \delta\mu}{\delta
(\delta\lambda- \kappa)}}$, (\ref{w2VI4}) is plainly true
for all $x \leq\sqrt{\frac{\kappa+ \delta\mu} {
\delta(\delta\lambda- \kappa)}}$.
On the other hand, we can use (\ref{tanh-alpha-ineqs2})
to calculate
\begin{eqnarray*}
f_2' (x) & =& \frac{\sqrt{2\delta}  [ \delta(\delta\lambda-
\kappa) x^2 - (\kappa+ \delta\mu)  ]} {\cosh\sqrt{2\delta}
x} \biggl[
\frac{\sqrt{2\delta} (\delta\lambda- \kappa) x}{\delta
(\delta\lambda- \kappa) x^2 - (\kappa+ \delta\mu)} - \cosh \sqrt{2\delta} x \biggr]
\\
& > &0 \qquad\mbox{for all } x \in\,\biggl ] \sqrt{\frac{\kappa+
\delta\mu} {\delta(\delta\lambda- \kappa)}} , \alpha\biggr [ ,
\end{eqnarray*}
and (\ref{w2VI2}) follows.\vadjust{\goodbreak}

Inequality (\ref{w2VI3}) is equivalent to
\[
\lambda- (\delta\lambda- \kappa) x^2 + \mu\leq0 \qquad \mbox{for all
} x \in\,\biggl ] \alpha, \frac{1}{2\lambda} \biggr[ \quad\Leftrightarrow\quad\alpha\geq\sqrt{
\frac{\lambda+ \mu} {
\delta\lambda- \kappa}} .
\]
In view of (\ref{tanh-alpha-ineqs2}), this is true if and only if
\[
\tanh\sqrt{\frac{2\delta(\lambda+ \mu)}{\delta\lambda-
\kappa}} < \sqrt{\frac{2\delta(\lambda+ \mu)}{\delta\lambda
- \kappa}}
\]
because $\sqrt{\frac{2\delta(\lambda+ \mu)}{\delta\lambda
- \kappa}} > \sqrt{\frac{\kappa+ \delta\mu} {\delta(\delta
\lambda- \kappa)}}  \Leftrightarrow \delta\lambda- \kappa
> 0$.
This inequality is indeed true because $\sqrt{\frac{2\delta
(\lambda+ \mu)}{\delta\lambda- \kappa}} > 1  \Leftrightarrow
 \delta\lambda+ \kappa+ 2 \delta\mu> 0$, and (\ref{w2VI3})
follows.
\end{pf*}

\begin{pf*}{Proof of Proposition~\ref{prop:caseI3}}
If we denote by $f_3 (\alpha)$ the right-hand side of (\ref{alpha3}),
then we can check that
%
\begin{equation}
\label{f3'} \qquad f_3' (\alpha) = - \frac{\sqrt{2\delta} \kappa [ \delta
(\delta
\lambda- \kappa) \alpha^2 - (\kappa+ \delta\mu)  ] + \delta
\sqrt{2\delta} (\delta\lambda- \kappa) (\delta- 2\kappa\alpha)
\alpha} { [ \delta(\delta\lambda- \kappa) \alpha^2 - (\kappa
+ \delta\mu)  ] ^2}
\end{equation}
{and}
\[
f_3'' (\alpha) = \frac{\delta\sqrt{2\delta} (\delta\lambda-
\kappa) (6 \kappa\alpha- \delta)} { [ \delta(\delta\lambda
- \kappa) \alpha^2 - (\kappa+ \delta\mu)  ] ^2} +
\frac{4 \delta^2 \sqrt{2\delta} (\delta\lambda- \kappa)^2
(\delta- 2 \kappa\alpha) \alpha} { [ \delta(\delta\lambda
- \kappa) \alpha^2 - (\kappa+ \delta\mu)  ] ^3} .
\]

If $\frac{\delta}{2\kappa} < \sqrt{\frac{\kappa+ \delta\mu} {
\delta(\delta\lambda- \kappa)}}$, then these calculations
imply that
\[
f_3' (\alpha) > 0 \quad\mbox{and} \quad
f_3'' (\alpha) > 0 \qquad \mbox{for all
} \alpha\in \biggl] \frac{\delta}{2\kappa} , \sqrt{\frac{\kappa+ \delta\mu} {\delta(\delta\lambda- \kappa)}} \biggr[ .
\]
Combining these inequalities with the observations that
\begin{eqnarray*}
f_3 (\alpha) &<& 0 \qquad\mbox{for all } \alpha\in \biggl[ 0,
\frac
{\delta} {
2\kappa} \biggr[ \,\cup\, \biggl] \sqrt{\frac{\kappa+ \delta\mu} {
\delta(\delta\lambda- \kappa)}} , \infty \biggr[ ,
\\[-3pt]
f_3 \biggl( \frac{\delta}{2\kappa} \biggr) &=& 0 \quad\mbox{and} \quad
\lim_{\alpha\uparrow\sqrt{\fracf{\kappa+ \delta\mu}
{\delta(\delta\lambda- \kappa)}}} f_3 (\alpha) = \infty
\end{eqnarray*}
and the fact that the restriction of $\tanh$ in $\R_+$ is
strictly concave, we can see that equation (\ref{alpha3})
has a unique solution $\alpha> 0$, which satisfies
(\ref{alpha3range1}).
In particular, we can see that
%
\begin{eqnarray}
&&\tanh\sqrt{2\delta} \alpha- \frac{\delta(\delta- 2\kappa
\alpha)} {\sqrt{2\delta}  [ \delta(\delta\lambda- \kappa)
\alpha^2 - (\kappa+ \delta\mu)  ]} \nonumber
\\[-9pt]
\\[-9pt]
&&\qquad \cases{ \displaystyle
> 0 , &\quad if $\displaystyle x \in \,]0, \alpha[ \,\cup \,\biggl] \sqrt{\frac{\kappa+ \delta\mu} {\delta(\delta\lambda
- \kappa)}} , \infty \biggr[$,\vspace*{3pt}
\cr
\displaystyle < 0 , &\quad if $\displaystyle x \in \,\biggl] \alpha, \sqrt{\frac{\kappa+ \delta\mu} {\delta(\delta
\lambda- \kappa)}} \biggr[ $.
} %
\label{tanh-alpha-ineqs3}
\nonumber
\end{eqnarray}

If $\sqrt{\frac{\kappa+ \delta\mu} {\delta(\delta\lambda
- \kappa)}} < \frac{\delta}{2\kappa}$, then (\ref{f3'}) implies
that
\[
f_3' (\alpha) < 0 \qquad\mbox{for all } \alpha\in \,\biggl]
\sqrt{\frac{\kappa+ \delta\mu} {\delta(\delta\lambda-
\kappa)}} , \frac{\delta}{2\kappa} \biggr[ .
\]
This inequality and the calculations
\begin{eqnarray*}
f_3 (\alpha) &<& 0 \qquad\mbox{for all } x \in \biggl[ 0, \sqrt{
\frac{\kappa+ \delta\mu} {\delta(\delta\lambda-
\kappa)}} \biggr[ \,\cup\,\biggl ] \frac{\delta} {2\kappa} , \infty \biggr[ ,
\\
\lim_{\alpha\downarrow\sqrt{\fracf{\kappa+ \delta\mu}
{\delta(\delta\lambda- \kappa)}}} f_3 (\alpha) &=& \infty \quad
\mbox{and} \quad f_3 \biggl( \frac{\delta} {2\kappa} \biggr) = 0 ,
\end{eqnarray*}
imply that equation (\ref{alpha3}) has a unique solution
$\alpha$ satisfying (\ref{alpha3range2}).
In particular, we can see that
\[
\frac{1}{2\lambda} < \alpha\quad\Leftrightarrow\quad \mbox{(\ref{cond32}) is
true} .
\]

We will show that the function $u$ satisfies all of the
requirements of Definition~\ref{u} if and only if we prove
that
%
\begin{equation}
\bigl\llvert u'(x)\bigr\rrvert \leq1 \qquad\mbox{for all } \llvert
x\rrvert \leq\alpha. \label{w3VI1}
\end{equation}
If the parameters are such that (\ref{alpha3range2}) is
true, then this inequality follows immediately from the
boundary conditions $u'(0) = 0$, $u_-'(\alpha) =1$
and the fact that $u$ is convex, which is true because
$A>0$.
If the parameters are such that (\ref{alpha3range1}) is
true, then \mbox{$A<0$}.
In this case, $u'''(x) = (2\delta)^{\fraca{3}{2}} A \sinh
\sqrt{2\delta} x < 0$ for all $x \in[0, \alpha[$, which
implies that $u''$ is strictly decreasing in $[0,\alpha[$.
Combining this observation with the fact that $u$ is
an even function, we can see that (\ref{w3VI1}) is true
if and only if $\lim_{x \uparrow\alpha} u'' (x) \geq0$,
which is equivalent to $\alpha\geq\sqrt{\frac{\mu} {
\delta\lambda- \kappa}}$.
In view of (\ref{tanh-alpha-ineqs3})\vspace*{-3pt} and the fact that
$\sqrt{\frac{\mu}{\delta\lambda- \kappa}} <
\sqrt{\frac{\kappa+ \delta\mu} {\delta(\delta\lambda
- \kappa)}}$, we\vspace*{2pt} can see that this indeed the case if
and only if (\ref{cond31}) is true.
\end{pf*}

\section{Proofs of results in Section~\texorpdfstring{\protect\ref{16161616}}{6}}\label{app2}
\begin{pf*}{Proof of Proposition~\ref{prop:caseII1}}
In view of (\ref{WSC21}), we will prove that $u$ satisfies
the conditions of Definition~\ref{u} if we show that
%
\begin{eqnarray}
u'(x) &\geq&-1 \qquad \mbox{for all } x \geq0 , \label{u21VI1}
\\
u(x) &\geq&- \lambda x^2 + \lambda\qquad \mbox{for all } x \geq
\alpha\label{u21VI2}
\end{eqnarray}
{and}
\begin{eqnarray}
\Lop u(x) + h(x) &\equiv&{\tfrac{1}{2}}u''(x) -
\delta u(x)
\nonumber
\\[-8pt]
\\[-8pt]
&\leq&0 \qquad \mbox{for all } x \in[0, \alpha] . \label{u21VI3}
\nonumber
\end{eqnarray}
Inequality (\ref{u21VI2}) follows immediately by the facts
that $u$ is $C^1$ at $\alpha$ and the restriction of $x
\mapsto u(x) + \lambda x^2 - \lambda$ in $[\alpha, \infty[$
is strictly convex.
Inequality (\ref{u21VI3}) is equivalent to $x^2 \leq
1 + \delta^{-1}$ for all $x \in[0, \alpha]$, which is true
because $\alpha< 1$.
Finally, inequality (\ref{u21VI1}) is true if and only if
$u'(\alpha) \geq-1$ because the restriction of $u'$ in
$[0,\infty[$ has a global minimum at $\alpha$.
Combining this observation with the identity $u' (\alpha)
= - 2\lambda\alpha$ and (\ref{Aalpha21}), we can see
that (\ref{u21VI1}) is satisfied if and only if (\ref{cond21})
true.
\end{pf*}

\begin{pf*}{Proof of Proposition~\ref{prop:caseII2}}
It is a matter of straightforward algebra to verify that $\alpha
> \frac{1}{2\lambda}$ if and only if the first inequality in
(\ref{cond22}) is true, which we assume in what follows.
Similarly, it is a matter of algebraic manipulations to show
that the constant on the left-hand side of (\ref{cond22}) is
strictly less than the constant on the right-hand side of
(\ref{cond22}).
Combining the inequality $\alpha> \frac{1}{2\lambda}$
with the strict concavity of the function $x \mapsto\lambda
x^2$, we can see that there exists $\beta\in[0,\alpha[$
such that the function $u$ defined by (\ref{u22}) is continuous
and $u(x) < \lambda x^2$ for all $x \in\, ]\beta,
\alpha[$ if and only if $\lambda\leq- \lambda\alpha^2 +
\alpha+ \lambda$, which is equivalent to the second
inequality in (\ref{cond22}).

We now assume that the problem data is such that
(\ref{cond22}) is true.
In view of the arguments above and (\ref{WSC22}), we
will prove that $u$ satisfies the requirements of Definition~\ref{u}
if we show that
%
\begin{eqnarray}
u'(x) &\geq&-1 \qquad \mbox{for all } x \in[0, \beta[ \, \cup
\, ]\alpha, \infty[ , \label{u22VI1}
\\
u(x) &\geq&- \lambda x^2 + \lambda \qquad \mbox{for all } x \geq
\alpha, \label{u22VI2}
\end{eqnarray}
{and}
\begin{eqnarray}
\Lop u(x) + h(x) &\equiv&{\tfrac{1}{2}}u''(x) -
\delta u(x)
\nonumber
\\[-8pt]
\\[-8pt]
&\leq&0 \qquad \mbox{for all } x \in[0, \beta] . \label{u22VI3}
\nonumber
\end{eqnarray}
The inequalities (\ref{u22VI1}) and (\ref{u22VI2}) follow
immediately by the facts that $u$ is $C^1$ at $\alpha$,
the restriction of $x \mapsto u(x) + \lambda x^2 - \lambda$
in $[\alpha, \infty[$ is strictly convex and $0 < \beta<
\frac{1}{2\lambda} < \alpha< 1$.
Finally, the inequality (\ref{u22VI3}) is equivalent to $x^2
\leq1 + \delta^{-1}$ for all $x \in[0, \beta]$, which is
plainly true because $\beta< 1$.
\end{pf*}

\begin{pf*}{Proof of Proposition~\ref{prop:caseII3}}
The inequality $u(x) < \lambda x^2$ for all $x \in[0,\alpha[$
that characterizes the region ${\mathcal S} _c = [-\alpha,
\alpha]$ is true if and only if $\lambda> - \lambda\alpha^2
+ \alpha+ \lambda$, which is equivalent to (\ref{cond23}).
Otherwise, the proof of this result is very similar to the
proof of Proposition~\ref{prop:caseII2}.
\end{pf*}
\end{appendix}



%

\printaddresses


\begin{thebibliography}{24}

\bibitem{BC}
%
\begin{bincollection}[mr]
\bauthor{\bsnm{Bather},~\bfnm{John}\binits{J.}} \AND
\bauthor{\bsnm{Chernoff},~\bfnm{Herman}\binits{H.}}
(\byear{1967}).
\btitle{Sequential decisions in the control of a spaceship}.
In \bbooktitle{Proc. {F}ifth {B}erkeley {S}ympos. {M}athematical
{S}tatistics and {P}robability ({B}erkeley, {C}alif., 1965/66), {V}ol.
{III}: {P}hysical {S}ciences}
\bpages{181--207}.
\bpublisher{Univ. California Press},
\blocation{Berkeley, CA.}
\bid{mr={0224218}}
\end{bincollection}
%
\bptok{imsref}%
\endbibitem

\bibitem{BH}
%
\begin{barticle}[auto:STB|2014/02/12|14:17:21]
\bauthor{\bsnm{Bayraktar},~\bfnm{E.}\binits{E.}} \AND
\bauthor{\bsnm{Huang},~\bfnm{Y.-J.}\binits{Y.-J.}}
(\byear{2012}).
\btitle{On the multi-dimensional controller and stopper games}.
\bjournal{SIAM J. Control Optim.}
\bvolume{51}
\bpages{1263--1297}.
\bid{mr={3036989}}
\end{barticle}
%
\bptok{imsref}%
\endbibitem

\bibitem{BY}
%
\begin{barticle}[mr]
\bauthor{\bsnm{Bayraktar},~\bfnm{Erhan}\binits{E.}} \AND
\bauthor{\bsnm{Young},~\bfnm{Virginia~R.}\binits{V.~R.}}
(\byear{2011}).
\btitle{Proving regularity of the minimal probability of ruin via a
game of stopping and control}.
\bjournal{Finance Stoch.}
\bvolume{15}
\bpages{785--818}.
\bid{doi={10.1007/s00780-011-0160-1}, issn={0949-2984}, mr={2863643}}
\end{barticle}
%
\bptok{imsref}%
\endbibitem

\bibitem{CH}
%
\begin{barticle}[mr]
\bauthor{\bsnm{Chiarolla},~\bfnm{Maria~B.}\binits{M.~B.}} \AND
\bauthor{\bsnm{Haussmann},~\bfnm{Ulrich~G.}\binits{U.~G.}}
(\byear{1998}).
\btitle{Optimal control of inflation: A central bank problem}.
\bjournal{SIAM J. Control Optim.}
\bvolume{36}
\bpages{1099--1132 (electronic)}.
\bid{doi={10.1137/S036301299630495X}, issn={0363-0129}, mr={1613921}}
\end{barticle}
%
\bptok{imsref}%
\endbibitem

\bibitem{DZ94}
%
\begin{barticle}[mr]
\bauthor{\bsnm{Davis},~\bfnm{M.~H.~A.}\binits{M.~H.~A.}} \AND
\bauthor{\bsnm{Zervos},~\bfnm{M.}\binits{M.}}
(\byear{1994}).
\btitle{A problem of singular stochastic control with discretionary stopping}.
\bjournal{Ann. Appl. Probab.}
\bvolume{4}
\bpages{226--240}.
\bid{issn={1050-5164}, mr={1258182}}
\end{barticle}
%
\bptok{imsref}%
\endbibitem

\bibitem{EMP}
%
\begin{barticle}[mr]
\bauthor{\bsnm{Ekeland},~\bfnm{Ivar}\binits{I.}},
\bauthor{\bsnm{Mbodji},~\bfnm{Oumar}\binits{O.}} \AND
\bauthor{\bsnm{Pirvu},~\bfnm{Traian~A.}\binits{T.~A.}}
(\byear{2012}).
\btitle{Time-consistent portfolio management}.
\bjournal{SIAM J. Financial Math.}
\bvolume{3}
\bpages{1--32}.
\bid{doi={10.1137/100810034}, issn={1945-497X}, mr={2968026}}
\end{barticle}
%
\bptok{imsref}%
\endbibitem

\bibitem{EC}
%
\begin{barticle}[auto:STB|2014/02/12|14:17:21]
\bauthor{\bsnm{El Karoui},~\bfnm{N.}\binits{N.}} \AND
\bauthor{\bsnm{Chaleyat-Maurel},~\bfnm{M.}\binits{M.}}
(\byear{1978}).
\btitle{Un probl\`eme de r\'eflexion et ses applications au temps
local at aux equations diff\'erentielles stochastiques sur $\R$, cas continu}.
\bjournal{Ast\'erisque}
\bvolume{52-53}
\bpages{117--144}.
\end{barticle}
%
\bptok{imsref}%
\endbibitem

\bibitem{H}
%
\begin{barticle}[mr]
\bauthor{\bsnm{Hamad{\`e}ne},~\bfnm{S.}\binits{S.}}
(\byear{2006}).
\btitle{Mixed zero-sum stochastic differential game and {A}merican
game options}.
\bjournal{SIAM J. Control Optim.}
\bvolume{45}
\bpages{496--518}.
\bid{doi={10.1137/S036301290444280X}, issn={0363-0129}, mr={2246087}}
\end{barticle}
%
\bptok{imsref}%
\endbibitem

\bibitem{HL}
%
\begin{barticle}[mr]
\bauthor{\bsnm{Hamad{\`e}ne},~\bfnm{S.}\binits{S.}} \AND
\bauthor{\bsnm{Lepeltier},~\bfnm{J.-P.}\binits{J.-P.}}
(\byear{2000}).
\btitle{Reflected {BSDE}s and mixed game problem}.
\bjournal{Stochastic Process. Appl.}
\bvolume{85}
\bpages{177--188}.
\bid{doi={10.1016/S0304-4149(99)00072-1}, issn={0304-4149}, mr={1731020}}
\end{barticle}
%
\bptok{imsref}%
\endbibitem

\bibitem{JJZ}
%
\begin{barticle}[mr]
\bauthor{\bsnm{Jack},~\bfnm{Andrew}\binits{A.}},
\bauthor{\bsnm{Johnson},~\bfnm{Timothy~C.}\binits{T.~C.}} \AND
\bauthor{\bsnm{Zervos},~\bfnm{Mihail}\binits{M.}}
(\byear{2008}).
\btitle{A singular control model with application to the goodwill problem}.
\bjournal{Stochastic Process. Appl.}
\bvolume{118}
\bpages{2098--2124}.
\bid{doi={10.1016/j.spa.2008.01.001}, issn={0304-4149}, mr={2462291}}
\end{barticle}
%
\bptok{imsref}%
\endbibitem

\bibitem{K}
%
\begin{barticle}[mr]
\bauthor{\bsnm{Karatzas},~\bfnm{Ioannis}\binits{I.}}
(\byear{1983}).
\btitle{A class of singular stochastic control problems}.
\bjournal{Adv. in Appl. Probab.}
\bvolume{15}
\bpages{225--254}.
\bid{doi={10.2307/1426435}, issn={0001-8678}, mr={0698818}}
\end{barticle}
%
\bptok{imsref}%
\endbibitem

\bibitem{KS}
%
\begin{barticle}[mr]
\bauthor{\bsnm{Karatzas},~\bfnm{Ioannis}\binits{I.}} \AND
\bauthor{\bsnm{Sudderth},~\bfnm{William~D.}\binits{W.~D.}}
(\byear{2001}).
\btitle{The controller-and-stopper game for a linear diffusion}.
\bjournal{Ann. Probab.}
\bvolume{29}
\bpages{1111--1127}.
\bid{doi={10.1214/aop/1015345598}, issn={0091-1798}, mr={1872738}}
\end{barticle}
%
\bptok{imsref}%
\endbibitem

\bibitem{KW}
%
\begin{barticle}[mr]
\bauthor{\bsnm{Karatzas},~\bfnm{I.}\binits{I.}} \AND
\bauthor{\bsnm{Wang},~\bfnm{H.}\binits{H.}}
(\byear{2000}).
\btitle{A barrier option of {A}merican type}.
\bjournal{Appl. Math. Optim.}
\bvolume{42}
\bpages{259--279}.
\bid{doi={10.1007/s002450010013}, issn={0095-4616}, mr={1795611}}
\end{barticle}
%
\bptok{imsref}%
\endbibitem

\bibitem{KZ2}
%
\begin{barticle}[mr]
\bauthor{\bsnm{Karatzas},~\bfnm{Ioannis}\binits{I.}} \AND
\bauthor{\bsnm{Zamfirescu},~\bfnm{Ingrid-Mona}\binits{I.-M.}}
(\byear{2008}).
\btitle{Martingale approach to stochastic differential games of
control and stopping}.
\bjournal{Ann. Probab.}
\bvolume{36}
\bpages{1495--1527}.
\bid{doi={10.1214/07-AOP367}, issn={0091-1798}, mr={2435857}}
\end{barticle}
%
\bptok{imsref}%
\endbibitem

\bibitem{LZ}
%
\begin{barticle}[mr]
\bauthor{\bsnm{Lamberton},~\bfnm{Damien}\binits{D.}} \AND
\bauthor{\bsnm{Zervos},~\bfnm{Mihail}\binits{M.}}
(\byear{2013}).
\btitle{On the optimal stopping of a one-dimensional diffusion}.
\bjournal{Electron. J. Probab.}
\bvolume{18}
\bpages{49}.
\bid{doi={10.1214/EJP.v18-2182}, issn={1083-6489}, mr={3035762}}
\end{barticle}
%
\bptok{imsref}%
\endbibitem

\bibitem{MS}
%
\begin{bincollection}[mr]
\bauthor{\bsnm{Maitra},~\bfnm{Ashok~P.}\binits{A.~P.}} \AND
\bauthor{\bsnm{Sudderth},~\bfnm{William~D.}\binits{W.~D.}}
(\byear{1996}).
\btitle{The gambler and the stopper}.
In \bbooktitle{Statistics, Probability and Game Theory}
(\beditor{T.~S. Ferguson, L.~S. Shapley and J.~B. MacQueen}, {eds.}).
\bseries{Institute of Mathematical Statistics Lecture
Notes---Monograph Series}
\bvolume{30}
\bpages{191--208}.
\bpublisher{IMS},
\blocation{Hayward, CA}.
\bid{doi={10.1214/lnms/1215453573}, mr={1481781}}
\end{bincollection}
%
\bptok{imsref}%
\endbibitem

\bibitem{MZ}
%
\begin{barticle}[auto:STB|2014/02/12|14:17:21]
\bauthor{\bsnm{Miller},~\bfnm{M.}\binits{M.}} \AND
\bauthor{\bsnm{Zhang},~\bfnm{L.}\binits{L.}}
(\byear{1996}).
\btitle{Optimal target zones: How an exchange rate mechanism can
improve upon discretion}.
\bjournal{J. Econom. Dynam. Control}
\bvolume{20}
\bpages{1641--1660}.
\end{barticle}
%
\bptok{imsref}%
\endbibitem

\bibitem{PS}
%
\begin{bbook}[mr]
\bauthor{\bsnm{Peskir},~\bfnm{Goran}\binits{G.}} \AND
\bauthor{\bsnm{Shiryaev},~\bfnm{Albert}\binits{A.}}
(\byear{2006}).
\btitle{Optimal Stopping and Free-Boundary Problems}.
\bpublisher{Birkh\"auser},
\blocation{Basel}.
\bid{mr={2256030}}
\end{bbook}
%
\bptok{imsref}%
\endbibitem

\bibitem{P}
%
\begin{bbook}[auto:STB|2014/02/12|14:17:21]
\bauthor{\bsnm{Protter},~\bfnm{P.}\binits{P.}}
(\byear{1990}).
\btitle{Stochastic Integration and Differential Equations}.
\bpublisher{Springer},
\blocation{Berlin}.
\bid{mr={1037262}}
\end{bbook}
%
\bptok{imsref}%
\endbibitem

\bibitem{RY}
%
\begin{bbook}[mr]
\bauthor{\bsnm{Revuz},~\bfnm{Daniel}\binits{D.}} \AND
\bauthor{\bsnm{Yor},~\bfnm{Marc}\binits{M.}}
(\byear{1994}).
\btitle{Continuous Martingales and {B}rownian Motion},
\bedition{2nd} ed.
\bseries{Grundlehren der Mathematischen Wissenschaften [Fundamental
Principles of Mathematical Sciences]}
\bvolume{293}.
\bpublisher{Springer},
\blocation{Berlin}.
\bid{mr={1303781}}
\end{bbook}
%
\bptok{imsref}%
\endbibitem

\bibitem{Sc}
%
\begin{barticle}[mr]
\bauthor{\bsnm{Schmidt},~\bfnm{W.}\binits{W.}}
(\byear{1989}).
\btitle{On stochastic differential equations with reflecting barriers}.
\bjournal{Math. Nachr.}
\bvolume{142}
\bpages{135--148}.
\bid{doi={10.1002/mana.19891420109}, issn={0025-584X}, mr={1017375}}
\end{barticle}
%
\bptok{imsref}%
\endbibitem

\bibitem{S}
%
\begin{barticle}[mr]
\bauthor{\bsnm{Sun},~\bfnm{Min}\binits{M.}}
(\byear{1987}).
\btitle{Singular control problems in bounded intervals}.
\bjournal{Stochastics}
\bvolume{21}
\bpages{303--344}.
\bid{doi={10.1080/17442508708833462}, issn={0090-9491}, mr={0905051}}
\end{barticle}
%
\bptok{imsref}%
\endbibitem

\bibitem{W}
%
\begin{barticle}[mr]
\bauthor{\bsnm{Weerasinghe},~\bfnm{Ananda}\binits{A.}}
(\byear{2006}).
\btitle{A controller and a stopper game with degenerate variance control}.
\bjournal{Electron. Commun. Probab.}
\bvolume{11}
\bpages{89--99 (electronic)}.
\bid{doi={10.1214/ECP.v11-1202}, issn={1083-589X}, mr={2231736}}
\end{barticle}
%
\bptok{imsref}%
\endbibitem

\bibitem{NW}
%
\begin{bincollection}[auto:STB|2014/02/12|14:17:21]
\bauthor{\bsnm{Williams},~\bfnm{N.}\binits{N.}}
(\byear{2008}).
\btitle{Robust control}.
In
\bbooktitle{The New Palgrave Dictionary of Economics}, \bedition{2nd} ed.
\bpublisher{Palgrave Macmillan},
\blocation{London}.
\end{bincollection}
%
\bptok{imsref}%
\endbibitem

\end{thebibliography}
\end{document}